\documentclass[mathpazo]{cicp}
\usepackage{tikz}
\usepackage{amssymb}
\usepackage{ragged2e}
\usepackage{indentfirst}
\usepackage{graphicx}
\usepackage{gensymb}
\usepackage{float}
\usepackage{subfigure}
\usepackage{booktabs}
\usepackage{array}
\usepackage{epstopdf}
\renewcommand{\raggedright}{\leftskip=0pt \rightskip=0pt plus 0cm}
\usepackage[colorlinks,
linkcolor=blue,
anchorcolor=blue,
citecolor=blue
]{hyperref}

\begin{document}
\title{A priori subcell limiting approach for the FR/CPR method on unstructured meshes}
\author{Guo-Quan Shi\affil{1} ,Huajun Zhu\affil{1,2}, Zhen-Guo Yan\affil{1}\comma\corrauth }
\address{\affilnum{1} State Key Laboratory of Aerodynamics, China Aerodynamics
	Research and Development Center, Mianyang, Sichuan 621000, China\\
	\affilnum{2} National University of Defense Technology, Changsha,
	Hunan 410073, China}
\emails{yanzhg@mail.ustc.edu.cn (Z.G. Yan)}
\begin{abstract}
	A priori subcell limiting approach is developed for high-order flux reconstruction/correction procedure via reconstruction (FR/CPR) on two-dimensional unstructured quadrilateral meshes. Firstly, a modified indicator based on modal energy coefficients is proposed to detect troubled cells. Then, troubled cells are decomposed into nonuniform subcells and each subcell has one solution point. A second-order finite difference shock-capturing scheme based on nonuniform nonlinear weighted (NNW) interpolation is constructed to calculate troubled cells while smooth cells are calculated by the CPR method. Numerical investigations show that the subcell limiting strategy on unstructured quadrilateral meshes is robust in shock-capturing. 
\end{abstract}
\maketitle

\section{Introduction}
	Although there are many industrial flow solvers based on second-order numerical methods, second-order methods have large dissipation and dispersion errors so that they are difficult to provide more accurate results for some complex flows than the high-order methods \cite{Wang2007}. Recently, high-order and high-resolution schemes have attracted increasing attention, especially, high order schemes on unstructured grids draw a lot of attention when the simulations involve complicated geometries, such as the k-exact finite volume (FV) methods, discontinuous Galerkin (DG) methods, spectral volume/spectral difference (SD/SV) methods \cite{Wang2007} and CPR method \cite{Huynh2007, Wang2009}, and some other methods\cite{Dumbser2010, Zhang2014}. 
\par	
	Among these high order methods, the CPR method is a differential method, which is both applicable to unstructured meshes and structured meshes and recovers a specific kind of DG, SD, as well as SV with appropriate choices of correction terms \cite{Huynh2007}. Moreover, the CPR method is more efficient than those methods which have expensive integration procedures, such as DG and SV\cite{Huynh2007, Yu201470}. The method is first proposed by Huynh in 2007 to solve the hyperbolic conservation law equations on the structured meshes, which was called the flux reconstruction (FR) method \cite{Huynh2007}. Wang and Gao extend FR to 2D triangular and hybrid meshes \cite{Wang2009}, proposing the lifting collocation penalty (LCP) method. Due to the close relationship between FR and LCP methods, the involved scholars named them the CPR method conformably \cite{Wang2011}. 
	The mathematical foundation \cite{Jameson2010,Castonguay2011,Vincent2010,Vincent2011,Williams2013,Williams2013-2} and simulations \cite{Wang2017,Jia2019,Li2013} using the CPR method have been widely studied. 
\par
	However, the CPR method is a high-order linear scheme, which probably generates spurious numerical oscillations near the discontinuities. Therefore, some shock-capturing strategies need to be used to suppress the numerical oscillations. 
	One approach is adding the artificial viscosity term, which is first proposed by von Neumann and Richtmyer \cite{VonNeumann1950}. In 2006, Persson and Peraire applied it for shock-capturing in the DG method \cite{Persson2006} and later Yu and Wang extended it to the FR method \cite{Yu2014}. Recently, Yu and Hesthaven study several artificial viscosity models within the DG framework \cite{Yu2020}. The artificial viscosity method has good accuracy retention characteristics, but it will increase the complexity of the equations. In addition, it is difficult to find universal parameters. 
\par
	There also exist some limiters to restrict the numerical oscillations. The multi-dimensional limiting process (MLP) originally in the FV framework has superior characteristics in terms of accuracy, robustness, and efficiency \cite{Kim2005}. Park and Kim extend MLP to the CPR method on unstructured meshes \cite{Park2016}. The MLP method can capture detailed flow structures in both continuous and discontinuous flow regions. However, the limited CPR loses its compactness because the stencil involves all the cells around the vertex.
	The WENO limiter of DG \cite{Qiu2005runge} also has a satisfactory shock-capturing performance, which gains many attractions \cite{Zhu2012RungeKuttaDG,Zhu2017RungeKuttaDG,Xu2009HierarchicalRF,Li2018,Li2020}. 
	Shu et al. also employ a WENO limiter on the CPR method, which can maintain high-order accuracy in smooth regions and control spurious numerical oscillations near discontinuities \cite{Du2015-1, Du2015}. 
\par 
	Another approach is the hybrid method, which adopts high-order methods in smooth regions to maintain compactness and high resolution, while FD, FV schemes, or lower-order FE schemes in troubled regions to provide robust shock-capturing abilities. 
	Cheng et al. present multi-domain hybrid RKDG/WENO-FD schemes based on domain decomposition \cite{Cheng2014AMH}. 
	In general, these FD, FV or FE schemes used for shock-capturing reduce the polynomial degree, so the resolution has to be preserved by h-refinement \cite{Dumbser2014, Dumbser2016, Guo2020, Sonntag2016, Hennemann2021}. 
	Dumbser et al. propose a DG/FV hybrid scheme, which divides a troubled cell into $2N+1$ uniform subcells ($N$ is the polynomial degree of DG) \cite{Dumbser2014, Dumbser2016}. Guo et al. propose a hybrid weighted compact nonlinear scheme and CPR (WCNS-CPR) scheme for simulating conservation laws\cite{Guo2020}, which adopts the WCNS scheme on the troubled cells. In these limiting strategies, the troubled cells are divided into uniform subcells. Since the solution points of these subcells do not coincide with the solution points of the DG or CPR approach, information exchange is necessary between the two kinds of solution points. Sonntag and Munz propose a subcell limiting strategy of DG \cite{Sonntag2016}, where the shock regions are treated by FV techniques. To keep the same integral mean value of the solution inside each DG cell, they choose the nodal values on Gauss points as the constant mean value of the subcells and Gauss weights as the length of the subcells. 
	Hennemann et al. discuss a provably entropy stable subcell shock capturing approach for DG, and the subcells are distributed nonuniformly \cite{Hennemann2021}. 
	These methods show good shock-capturing performance and can combine the advantages of different schemes and maintain a high resolution.
\par
	In this paper, we investigate a subcell limiting strategy of the high-order CPR method for two-dimensional unstructured quadrilateral meshes, it combines the shock-capturing ability of the second-order scheme and the advantages of the high computational efficiency of CPR. Similar to \cite{Sonntag2016}, we divide troubled cells into nonuniformly spaced subcells. We discuss the performance of the second-order shock-capturing scheme as a subcell shock-capturing scheme on unstructured meshes, which is originally applied on structured meshes in our recent works \cite{Zhu2021Shock}. In addition, we discuss and improve a troubled cell indicator for the subcell limiting strategy. 
\par
	A brief review of the CPR method is presented in Section \ref{CPR review}.
	Based on solution points of CPR, we develop a shock-capturing scheme based on second-order nonuniform nonlinear weighted interpolation in Section \ref{shock-capturing scheme}.
	The subcell limiting strategy is given in Section \ref{subcell limiting}, including and the strategy of the discontinuous detection and scheme switching. Several numerical tests about the performance of the proposed subcell limiting strategy are presented in Section \ref{numerical results}. Finally, concluding remarks are given in Section \ref{concluding}.
\section{Brief review of CPR}\label{CPR review}
	Consider the two-dimensional (2D) conservation law :
	\begin{equation}
		\frac{\partial U}{\partial t}+\frac{\partial F(U)}{\partial x}+\frac{\partial G(U)}{\partial y} =0, \label{2D-con-law}
	\end{equation}
	where $U$ is the conservative variable vector, $F$ and $G$ are the flux vectors. 
	For the 2D Euler equations, Eq.~\eqref{2D-con-law} becomes
	\begin{equation}\label{2D Euler}
	\frac{\partial}{\partial t}\left[\begin{array}{c}
	\rho \\
	\rho u \\
	\rho v \\
	E
	\end{array}\right]+\frac{\partial}{\partial x}\left[\begin{array}{c}
	\rho u \\
	\rho u^{2}+p \\
	\rho u v \\
	u(E+p)
	\end{array}\right]+\frac{\partial}{\partial y}\left[\begin{array}{c}
	\rho v \\
	\rho u v \\
	\rho v^{2}+p \\
	v(E+p)
	\end{array}\right]=0,
	\end{equation}
	with the equation of state 
	\begin{equation}
	p=(\gamma-1)\left[E-\frac{\rho}{2}\left(u^{2}+v^{2}\right)\right],
	\end{equation}
	where $\rho$ is the density, $u$ and $v$ are the velocities, $E$ is the total energy,	and $p$ is the pressure, and $\gamma$ is the ratio of specific heats with $\gamma=1.4$ for ideal gas.
\par
\begin{figure}[htbp]\centering
		\includegraphics[scale=0.6]{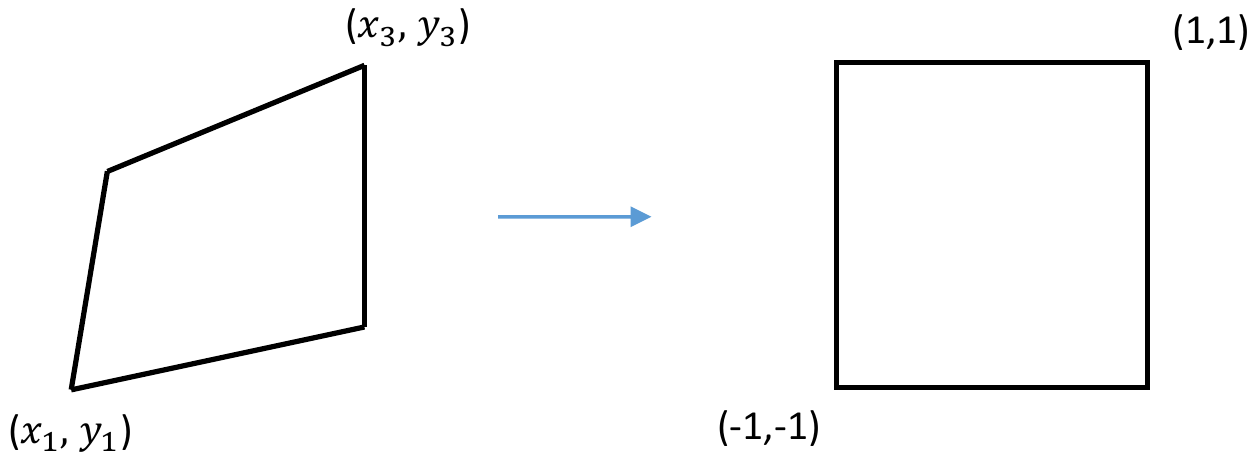}
		\caption{Transformation of a general straight quadrilateral cell to a standard cell.}\label{PhyCellToComCell}		
	\end{figure}

	The FR/CPR approach was first proposed by Huynh \cite{Huynh2007} and then extended to different mesh types by Wang and Gao \cite{Wang2009}. We give a brief review of CPR on unstructured quadrilateral meshes.  
	The $j$-th cell $E_j$ in the physical space is transformed into the standard cell $I_j=[-1,1]\times[-1,1]$, as shown in Figure \ref{PhyCellToComCell}. 
	The transformation can be written as
	\begin{equation}
		\binom{x}{y}=\sum_{i=1}^{K_v}M_{i}(\xi,\eta)\binom{x_j}{y_j},	\label{transform}
	\end{equation}
	where $K_v$ is the number of points used to define the physical cell, $M_{i}(\xi,\eta)$ are shape functions, and $(x_j,y_j)$ are the physical coordinates of those points. If the physical cell is a straight quadrilateral cell then $K_v$ equals 4 and $M_{i}$  follows
	\begin{equation}
		M_i=\frac{1}{4}(1+{\xi}_i\xi)(1+{\eta}_i\eta), i=1,2,3,4, \label{shape func}
	\end{equation}
	where $({\xi}_i, {\eta}_i)$ are the local coordinates in the standard element.\par
	The Jacobian matrix of the transformation is
	\begin{equation}
		J=\frac{\partial(x,y)}{\partial(\xi,\eta)}
		=\begin{bmatrix}
		x_{\xi} & x_{\eta}\\
		y_{\xi} & y_{\eta}
		\end{bmatrix},
	\end{equation}
	with $|J|={x}_{\xi}{y}_{\eta}-{x}_{\eta}{y}_{\xi}$.
	The governing equations \eqref{2D-con-law} in the physical domain are then transformed into the computational domain as:
	\begin{equation}
	\frac{\partial \tilde{U}}{\partial t} +\frac{\partial \tilde{F}}{\partial\xi } +\frac{\partial \tilde{G}}{\partial\eta}  =0, \label{Jacobi}
	\end{equation}
	where
	\begin{equation}		
		\tilde{U}=|J|\cdot U,~\tilde{F}=|J|({\xi}_{x}F+{\eta}_{y}G),~\tilde{G}=|J|({\eta}_{x}F+{\xi}_{y}G).
	\end{equation}
\par
	Based on the $U$ of the current time step, the fluxes on the solution points, $\tilde{F}_{j;k,l}$ and $\tilde{G}_{j;k,l}$, are first calculated. Then,
	flux polynomials are constructed by the tensor product of 1D Lagrange polynomials,
	\begin{equation}\label{F-Intepolate}
	\tilde{F}_j(\xi,\eta)=\sum_{k=1}^{N+1}\sum_{l=1}^{N+1}L_k(\xi)L_l(\eta)\tilde{F}_{j;k,l},
	\end{equation}
	\begin{equation}\label{G-Intepolate}
	\tilde{G}_j(\xi,\eta)=\sum_{k=1}^{N+1}\sum_{l=1}^{N+1}L_k(\xi)L_l(\eta)\tilde{G}_{j;k,l},
	\end{equation}
	where $N$ is the degree of the polynomials, the subscript $j$ is an index of element number, $L_k(\xi)$ and $L_l(\eta)$ are 1D Lagrange polynomials in the $\xi$ and $\eta$ directions, respectively. In our procedure, primary variables on solution points are calculated first. And then interpolations are performed based on the primary variables to get variables on element interfaces, which have a similar form as equation(\ref{F-Intepolate}) and (\ref{G-Intepolate}). 
\par
	For each fixed $\xi$ (or $\eta$), the interpolation along the $\xi$-direction (or $\eta$-direction) reduces to a 1D problem. To maintain continuity of the fluxes on element interfaces, a continuous flux function is constructed as
	\begin{equation}
		\tilde{F}_{j}^{con}(\xi,\eta_{l})=\tilde{F}_{j}(\xi,\eta_{l})+[\hat{F}_{j}(-1,\eta_{l})-\tilde{F}_{j}(-1,\eta_{l})]g_{L}(\xi)+[\hat{F}_{j}(1,\eta_{l})-\tilde{F}_{j}(1,\eta_{l})]g_{R}(\xi),  \label{con flux poly}
	\end{equation}
	where $\hat{F}(\pm 1,\eta_{l})=\hat{F}(U_{j;l},U_{j+;l},\vec{n})$ are the Riemann fluxes at the cell interfaces and $U_{j+}$ denote the interface values from the neighbor cells.  
	$g_{L}(\xi)$ and $g_{R}(\xi)$ are the correction functions, which are both degree $N+1$ polynomials satisfying :
	\begin{equation}	
		g_{L}(-1)=1,~g_{L}(1)=0,\quad
		g_{R}(-1)=0,~g_{R}(1)=1.	\label{cor func}
	\end{equation}
	$\tilde{G}_{j}^{con}(\xi_{k},\eta)$ can be constructed similarly in the $\eta$ direction.
	Therefore, the semi-discrete form of the Eq.~\eqref{Jacobi} is expressed as 
	\begin{equation}
	\begin{aligned}
		\frac{\partial \tilde{U}_j({\xi}_{k},{\eta}_{l})}{\partial t} = &
		-\frac{\partial \tilde{F}_{j}({\xi}_{k},{\eta}_{l})}{\partial \xi}
		-\frac{\partial \tilde{G}_{j}({\xi}_{k},{\eta}_{l})}{\partial \eta} \\
		&-((\hat{F}(-1,{\eta}_{l})-\tilde{F}(-1,{\eta}_{l})){g'}_{L}(\xi_{k})
		+(\hat{F}(1,{\eta}_{l})-\tilde{F}(1,{\eta}_{l})){g'}_{R}(\xi_{k}))\\
		&-((\hat{G}({\xi}_{k},-1)-\tilde{G}({\xi}_{k},-1)){g'}_{L}(\eta_{l})
		+(\hat{G}({\xi}_{k},1)-\tilde{G}({\xi}_{k},1)){g'}_{R}(\eta_{l}))
	\end{aligned}
	\end{equation}
	in which the expression of ${\partial \tilde{F}_j({\xi}_{k},{\eta}_{l})}/{\partial \xi}$ and ${\partial \tilde{G}_j({\xi}_{k},{\eta}_{l})}/{\partial \eta}$ can be obtained with the polynomial expression of $\tilde{F}$ and $\tilde{G}$ in Eqs.~\eqref{F-Intepolate} and~\eqref{G-Intepolate}. Combined with the explicit third-order TVD Runge-Kutta method \cite{Shu1988}, the solution $\tilde{U}_j({\xi}_{k}, {\eta}_{l})$ can be updated. We use the local Lax-Friedrichs flux as Riemann flux at the cell interfaces, and $g_{DG}$ function in \cite{Huynh2007} as the correction function. 
\section{Second order scheme based on NNW interpolation}\label{shock-capturing scheme}
	Recently, we developed shock-capturing schemes based on NNW interpolation, named compact NNW (CNNW), and applied them as subcell limiters for fifth-order CPR on structured meshes. Numerical results have shown good properties of the hybrid scheme \cite{Zhu2021Shock}, but only on structured meshes. In this paper, we generalize the second-order CNNW (CNNW2) to unstructured meshes and use them as a subcell limiting scheme for CPR.
	The length of each of subcells can be denoted by $\xi_{fp_{i+1}}-\xi_{fp_{i}}=\omega_{i},~\xi_{fp_{i}}=-1$, where $\omega_{i}$ are the weights of the Gauss integration, which is same as doing in \cite{Sonntag2016}. Here we study some technical details of subcell limiting on the unstructured quadrilateral meshes.
	\begin{figure}[htbp]\centering
		\includegraphics[scale=0.7]{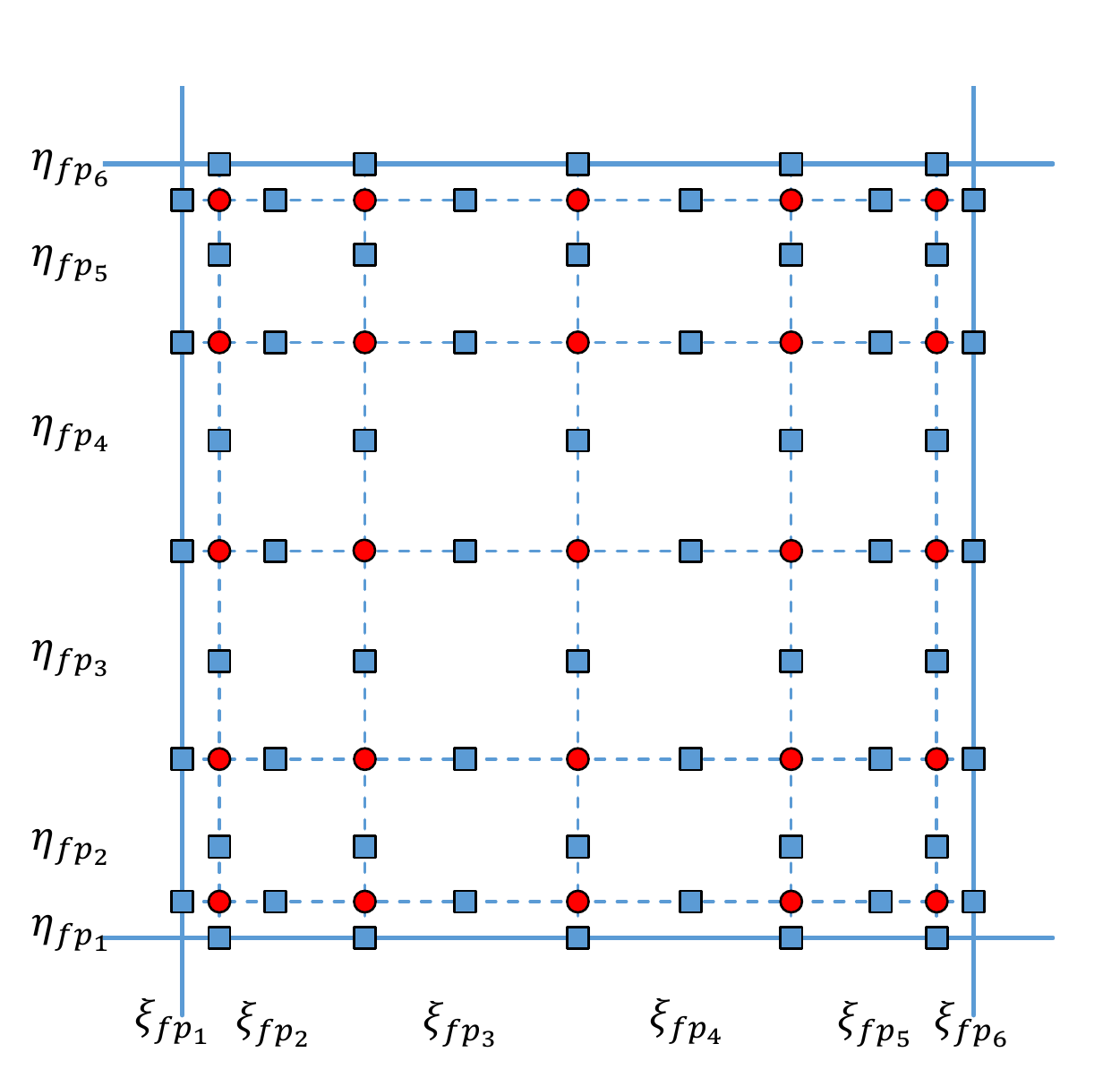}
		\caption{CPR reference cell is split into subcells --- with solution points marked in red dot and flux points in blue square.}\label{SPs and FPs}		
	\end{figure}
	\subsection{Nonuniform nonlinear weighted interpolation}\label{NNW I}
\begin{figure}[htbp]\centering
		\includegraphics[scale=0.8]{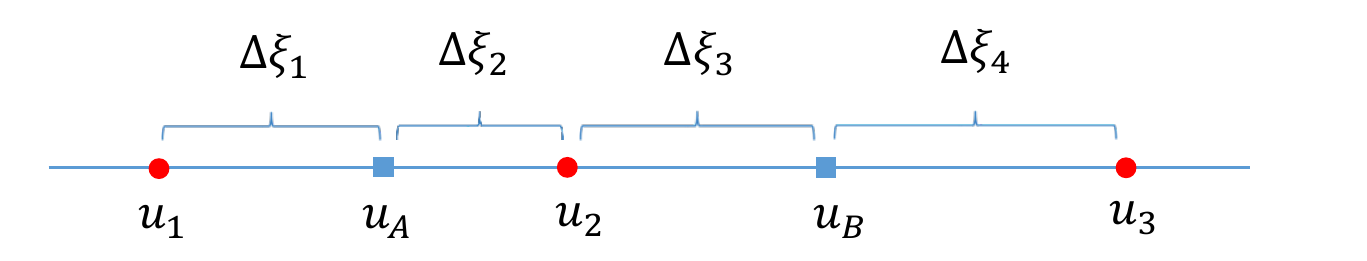}
		\caption{The stencil of nonuniform nonlinear interpolation for $u_A^{(1)}$ and $u_B^{(1)}$. $u_i,~i=1,2,3$ are the nodal solution values in the subcells. $u_A$ and $u_B$ are the subcell interfaces. $\Delta {\xi}_i,~i=1,2,3 $ are the distance between solution points and flux points.}\label{nonuniform cell}		
	\end{figure}

	We first discuss the one-dimensional case. Assume that there are $N_c$ uniform cells and we divide each cell into $K$ ($K$ is the number of solution points in a cell) subcells using the weights of the Gauss integration as subcell length. Consider a stencil with three nonuniformly distributed solution points $u_i,~i=1,2,3$, as shown in Figure \ref{nonuniform cell}.
	The values at the subcell interfaces $u_A$ and $u_B$ can be obtained by the following procedure :\\ 
	\\
	\begin{enumerate}
		\item 
		Calculate $u_A^{(1)}$ and $u_B^{(1)}$ by the inverse distance weighted interpolation.
			\begin{equation}
				\begin{aligned}
				\omega_{1} &=\frac{\left(1 / \Delta \xi_{1}\right)}{\left(1 / \Delta \xi_{1}\right)+\left(1 / \Delta \xi_{2}\right)}, \quad \omega_{2}=\frac{\left(1 / \Delta \xi_{2}\right)}{\left(1 / \Delta \xi_{1}\right)+\left(1 / \Delta \xi_{2}\right)}, \quad u_{A}^{(1)}=\omega_{1} u_{1}+\omega_{2} u_{2} , \\
				\omega_{3} &=\frac{\left(1 / \Delta \xi_{3}\right)}{\left(1 / \Delta \xi_{3}\right)+\left(1 / \Delta \xi_{4}\right)}, \quad \omega_{4}=\frac{\left(1 / \Delta \xi_{4}\right)}{\left(1 / \Delta \xi_{3}\right)+\left(1 / \Delta \xi_{4}\right)^{\prime}}, \quad u_{B}^{(1)}=\omega_{3} u_{2}+\omega_{4} u_{3},
				\end{aligned}\label{NNW2FistLayerInterplation}
			\end{equation}
			where ${\omega}_i$ are the interpolation weights based on the inverse distance between the solution points and the interpolation points A and B.\\
		\item 
		Calculate the gradient $\frac{\partial u}{\partial \xi}$ based on $u_A^{(1)}$, $u_2$ and $u_B^{(1)}$. 
		\begin{equation}	
				\frac{\partial u}{\partial \xi}  = \omega_{5}\left(\frac{\partial u}{\partial\xi}\right)^{(1)}+\omega_{6}\left(\frac{\partial u}{\partial \xi}\right)^{(2)}
				 = \omega_{5}\frac{u_{2}-u_{A}^{(1)}}{\Delta\xi_{2}}+\omega_{6}\frac{u_{B}^{(1)}-u_{2}}{\Delta\xi_{3}},	
			\label{eq:NNWI-derivative}
		\end{equation}
		in which 
		\begin{equation}
			\omega_{5}=\frac{\left(1/\Delta\xi_{2}\right)}{\left(1/\Delta\xi_{2}\right)+\left(1/\Delta\xi_{3}\right)},\quad\omega_{6}=\frac{\left(1/\Delta\xi_{3}\right)}{\left(1/\Delta\xi_{2}\right)+\left(1/\Delta\xi_{3}\right)}.
		\end{equation}
		In the derivation of Eq.~\eqref{eq:NNWI-derivative}, the inverse distance weighting method is adopted again to approximate the gradient based on approximations from both sides of the solution point.
		\item 
		Recompute $u_A$ and $u_B$ based on $u_2$ and the gradient ${\partial u}/{\partial \xi}$.
		\begin{equation}
		u_{A}^{(2)}=u_{2}-\frac{\partial u}{\partial \xi} \Delta \xi_{2},~u_{B}^{(2)}=u_{2}+\frac{\partial u}{\partial \xi} \Delta \xi_{3}.
		\end{equation}
		\item 
		Limit the gradient to control numerical oscillations. $u_A^R$ and $u_B^L$ are obtained by linear reconstruction with a limiter.
		\begin{equation}
			\begin{array}{c}
			u_{A}^{R}=u_{2}-\phi \frac{\partial u}{\partial \xi} \Delta \xi_{2},\qquad u_{B}^{L}=u_{2}+\phi \frac{\partial u}{\partial \xi} \Delta \xi_{3},\qquad
			\phi=min \left(lim(u_{A}^{(2)}), lim(u_{B}^{(2)})\right)
			\end{array}
		\end{equation}
		where the limiter \cite{Barth1989} is 
		\begin{equation}
			lim(u)=\left\{\begin{aligned}
			min \left(1, \frac{M-u_{2}}{u-u_{2}}\right), &if~u>u_{2}, \\
			min \left(1, \frac{m-u_{2}}{u-u_{2}}\right), &if~u<u_{2}, \\
			1, \qquad\qquad\qquad\quad &if u=u_{2}.
			\end{aligned}\right.
		\end{equation}
		with $m=min\left(u_{1},u_{2},u_{3}\right)$ and $M=max \left(u_{1},u_{2},u_{3}\right)$.
	\end{enumerate}
\par
	With the NNW interpolation, we can get the left and right values on the interfaces of the subcells. With the Riemann flux $\hat{F}_{j,fp_{k}}=\hat{F}(u_{j,fp_k},u_{j+,fp_k}, \vec{n})$, we can approximate the spatial derivative using
	\begin{equation}
	\begin{aligned}
		\frac{\partial F}{\partial \xi}= \frac{\hat{F}_{j,fp_{k+1}}-\hat{F}_{j,fp_{k}}}{\xi_{f p_{k+1}}-\xi_{f p_{k}}},
	\end{aligned}
	\end{equation}
	which is adopted to update the solutions combined with the explicit third-order TVD Runge-Kutta method. For two-dimensional case, the flux derivatives $\frac{\partial F}{\partial \xi},\frac{\partial G}{\partial \eta}$ can be discretized similarly.  
	\subsection{Extension to unstructured quadrilateral meshes}
\begin{figure}[htbp]\centering
		\includegraphics[scale=0.25]{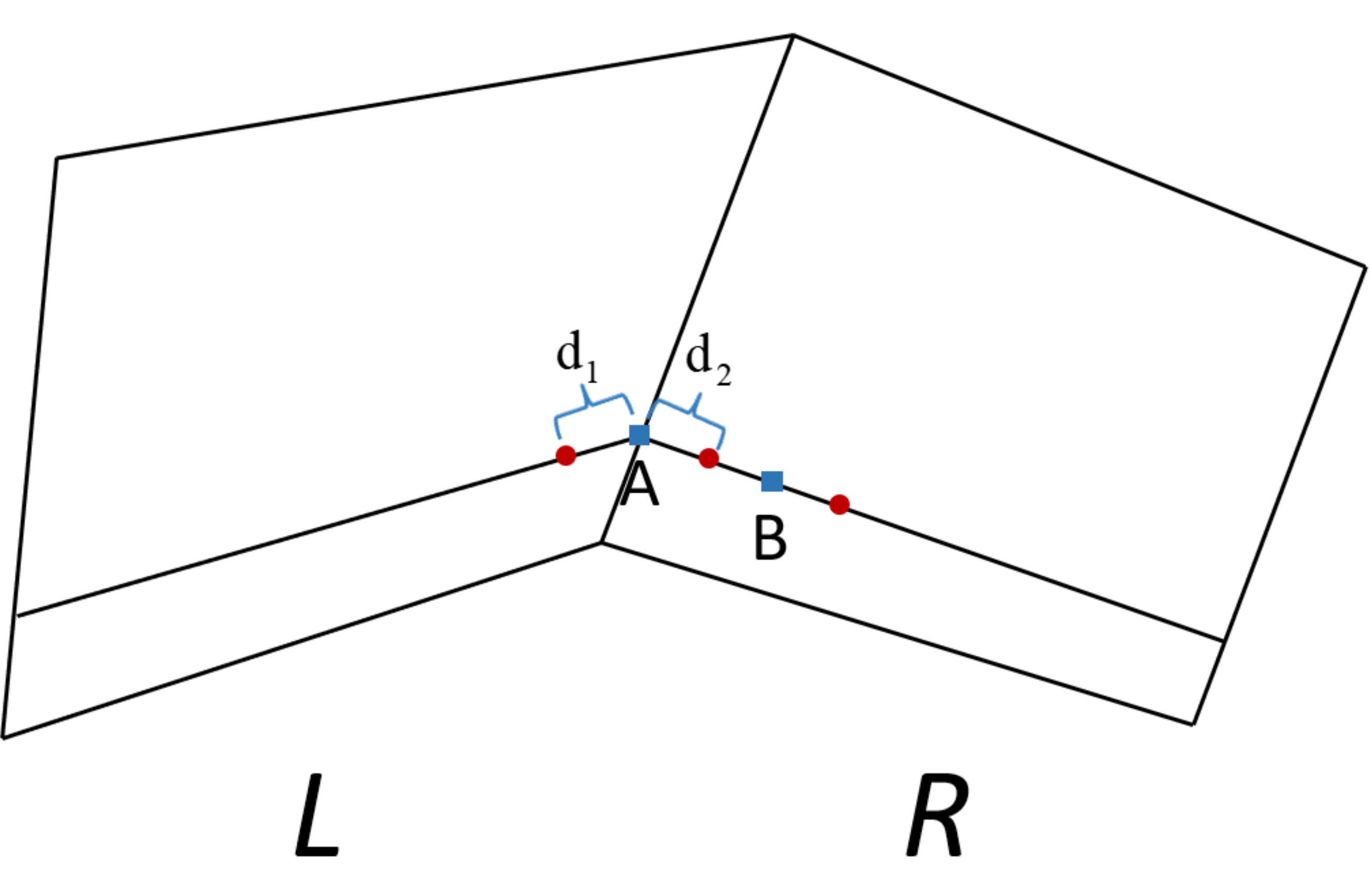}
		\caption{The distribution of nodal values of nonuniform nonlinear interpolation stencil in phycial domain. Solution points are marked in red dot and
			flux points in blue square.}\label{NNW Interpolate correction}		
	\end{figure}

	The interpolation approach presented in Subsection \ref{NNW I} offers an effective discontinuity regulation method based on slop limiting in the computational space, which however relies on a smooth grid transformation. Therefore, it needs modification for application on unstructured meshes, the grid transformation of which is usually not smooth across elements.
	The subcells are very similar to multi-block structured meshes, in which the subcells within one block are structured while relations between the blocks are unstructured.   
	In designing the subcell numerical scheme, we try to utilize the structured information within one subcell block and make modifications only for the cells next to block interfaces to adjust to unstructured meshes.
	Take the left side of the cell in Fig.~\ref{NNW Interpolate correction} as an example, computing $u_A^{(1)}$ needs data from the left neighbor cell. The distance of $d_1$ and $d_2$, which are always the same in computational space, could be quite different in physical space on unstructured meshes.
	Thus, instead of interpolating in computational space, we consider interpolation in physical space to obtain $u_{A}^{(1)}$, using
	\begin{equation}
	\omega_{1}=\frac{\left(1 / \Delta d_{1}\right)}{\left(1 / \Delta d_{1}\right)+\left(1 / \Delta d_{2}\right)}, \quad \omega_{2}=\frac{\left(1 / \Delta d_{2}\right)}{\left(1 / \Delta d_{1}\right)+\left(1 / \Delta d_{2}\right)}, \qquad u_{A}^{(1)}=\omega_{1} u_{1}+\omega_{2} u_{2},
	\end{equation}
	where $\Delta d_i,~i=1,2$ are the physical distance between the solution point and the flux point in the adjacent cell and the current cell, as shown in Fig.~\ref{NNW Interpolate correction}. This only changes the interpolation of Eqs.~\eqref{NNW2FistLayerInterplation} for subcells next to the element interfaces with all the other operations still performed in the computational space. 	
\section{Subcell limiting of CPR based on CNNW2}\label{subcell limiting}
	\subsection{Troubled cell detection}
	The subcell limiting strategy adopts the second-order scheme with shock-capturing ability in troubled regions and the CPR method with high computational efficiency in smooth regions. Detecting the troubled cells accurately is very important for a hybrid scheme to take full advantage of the hybridization.
	Recently, a shock indicator using modal energy of solution polynomials is developed. 
	As mentioned in \cite{Persson2006}, the modal coefficient will decay rapidly in smooth regions while slower in non-smooth regions. Therefore, the shock indicator compares the highest mode with a threshold value which is defined to decide whether a cell contains discontinuities. Hennemann et al. improved this method by using the highest and the second-highest mode to avoid odd/even effects \cite{Hennemann2021}. The definition of the modal energy of a 1D polynomial in \cite{Hennemann2021} is
	\begin{equation}
	\langle\epsilon, \epsilon\rangle_{L^{2}}=\left\langle\sum_{j=0}^{N} m_{j} \tilde{L}_{j}, \sum_{j=0}^{N} m_{j} \tilde{L}_{j}\right\rangle_{L^{2}}=\sum_{i, j=0}^{N} m_{i} m_{j}\left\langle\tilde{L}_{i}, \tilde{L}_{j}\right\rangle_{L^{2}}=\sum_{j=0}^{N} m_{j}^{2},
	\end{equation}
	where $\left\{m_{j}\right\}_{j=0}^{N}$ are the modal coefficients, $\tilde{L}_{i},i=0,1,...,N$ are the Legendre basis functions. To obtain $\left\{m_{j}\right\}_{j=0}^{N}$, we need transform Lagrange polynomial to Legendre polynomial with transformation matrix $\mathbb{K}$. We have
	\begin{subequations}\label{Lagrange Legendre}
		\begin{align}
			&\sum_{k=1}^{N+1} u_{k} L_{k}(\xi)=\sum_{k=0}^{N}A_k(u_1,...,u_{N+1},\xi_1,...,\xi_{N+1}){\xi}^k=A\Xi=\mathbb{A}\mathbb{U}\Xi,
			\\
			&\sum_{k=0}^{N} m_{k} \tilde{L}_{k}(\xi)=\sum_{k=0}^{N}M_k(m_0,...,m_{N}){\xi}^k=M\Xi=\mathbb{B}\mathbb{M}\Xi,
		\end{align}
	\end{subequations}
	where $\left\{u_{k}\right\}_{k=1}^{N+1}$ are the nodal values at the solution points, $L_{k}(\xi)$ are Lagrange basis functions, and $\mathbb{A}$ and $\mathbb{B}$ are both (N+1)-th order square matrices whose elements are related to the local coordinates of the solution points while $\mathbb{U}=(u_1,u_2,...,u_{N+1})^T$,  $\mathbb{M}=(m_0,m_1,...,m_{N})^T$ and $\Xi=(\xi^0,\xi^1,...,\xi^{N})$. 
	Thus,
	\begin{equation}\label{transform K}
	\mathbb{M}=\mathbb{B}^{-1}\mathbb{A}\mathbb{U}=\mathbb{K}\mathbb{U},
	\end{equation}
	and the matrix $\mathbb{K}$ is a constant matrix which can be solved in advance.
\par
	The percentage of the highest energy mode to the total energy of the polynomial is estimated as follows:
	\begin{equation}
	\mathbb{E}=\max \left(\frac{m_{N}^{2}}{\sum_{j=0}^{N} m_{j}^{2}}, \frac{m_{N-1}^{2}}{\sum_{j=0}^{N-1} m_{j}^{2}}\right).
	\end{equation}
	And a threshold value $\mathbb{T}=\mathbb{T}(N)$ is defined as
	\begin{equation}
	\mathbb{T}(N)=a \cdot 10^{-c(N+1)^{\frac{1}{4}}},
	\end{equation}
	where $N$ is the order of the polynomial, and a=0.5 and c=1.8. 
	A cell is marked as a troubled cell if $\mathbb{E}>\mathbb{T}(N)$. 
\par
\begin{figure}[H]	
		\centering
		\subfigbottomskip=2pt
		\subfigcapskip=-5pt
		\subfigure[]{		\label{ME ori}
			\includegraphics[width=0.45\linewidth]{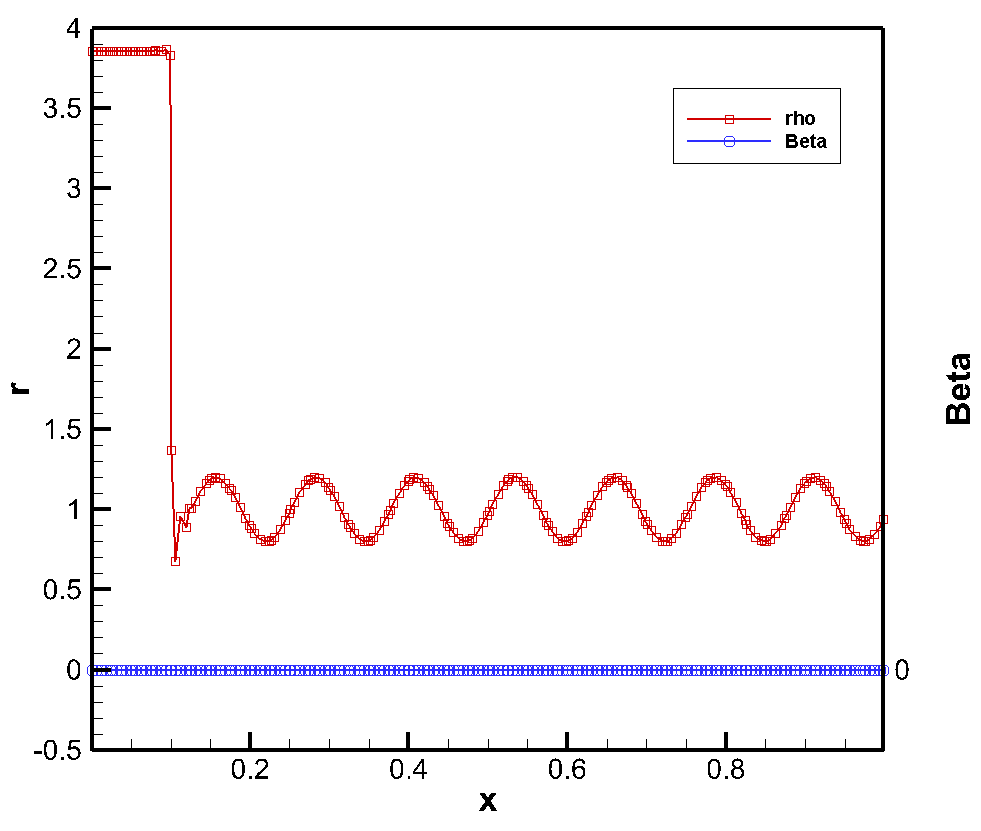}}
		\subfigure[]{
			\label{ME cor}
			\includegraphics[width=0.45\linewidth]{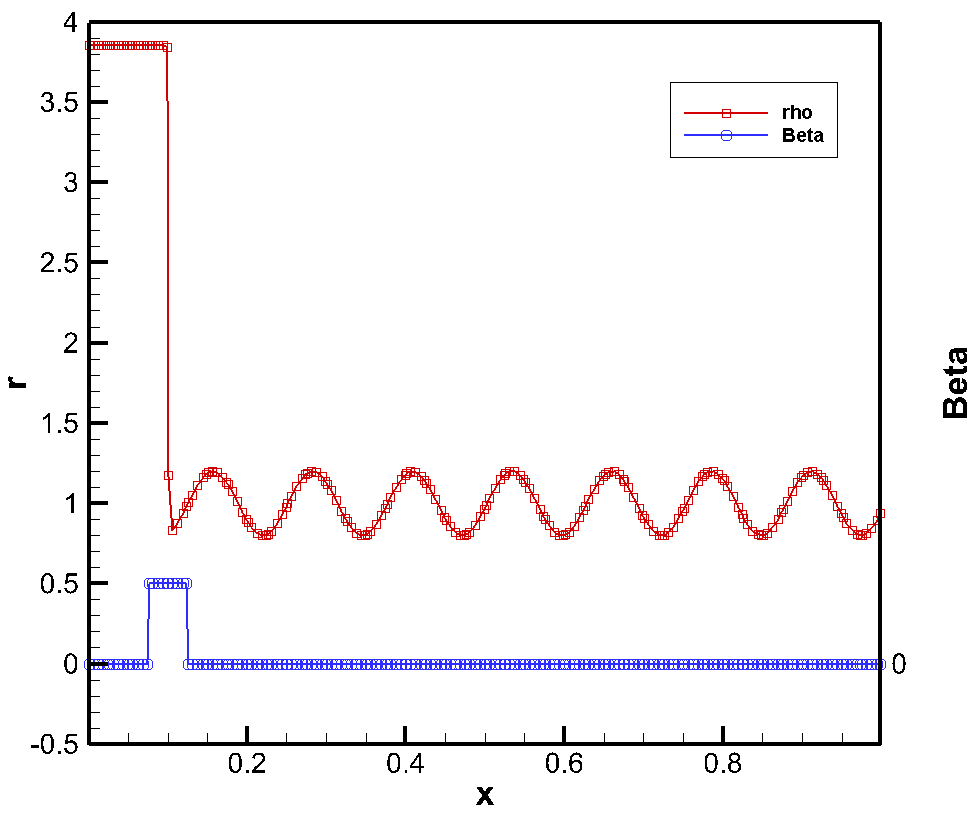}}	
		\caption{Results calculated by the hybrid CPR-CNNW2 scheme using modal decay indicator. (a) Original indicator (b) Improved indicator.}
	\end{figure}

 
	The shock indicator described above will fail when a strong discontinuity appears in the interface of two cells. We can see that the indicator fails to detect the shock at x = 0.1 from Figure \ref{ME ori}. 
	Thus, we consider the points of the cell interfaces together with the solution points inside the cell to construct a new polynomial. And using the Roe average as the interface values. The Roe average can be obtained by the two solution point values near the interface directly, 
	\begin{equation}
	\tilde{\rho}=\sqrt{\rho_L \rho_R}, \quad
	\tilde{f}=\frac{f_L\sqrt{\rho_L}+f_R\sqrt{\rho_R}}{\sqrt{\rho_L}+\sqrt{\rho_R} },\quad
	\tilde{p}=({\tilde{h}}^2-\frac{1}{2}\tilde{\rho}{\tilde{V}^2})\frac{(\gamma-1)}{\gamma},
	\end{equation}
	where $f=u,v,h$ and $\tilde{V}=\sqrt{{\tilde{u}}^2+{\tilde{v}}^2}$.
	The polynomial order will be increased to $N+2$ order, 
	\begin{equation}
		\sum_{k=1}^{N+3} u_{k} L_{k}(\xi)=\sum_{k=0}^{N+2} m_{k} \tilde{L}_{k}(\xi),
	\end{equation}
	and the other processes are the same to equation(\ref{Lagrange Legendre}) and equation(\ref{transform K}). We can see the significant improvement from Figure \ref{ME ori} to Figure \ref{ME cor}, which shows the first time stage of Runge Kutta computational results at the first time step. 
	\subsection{Calculation of Riemann flux at cell interfaces}
\begin{figure}[htbp]\centering	
			\includegraphics[scale=0.4]{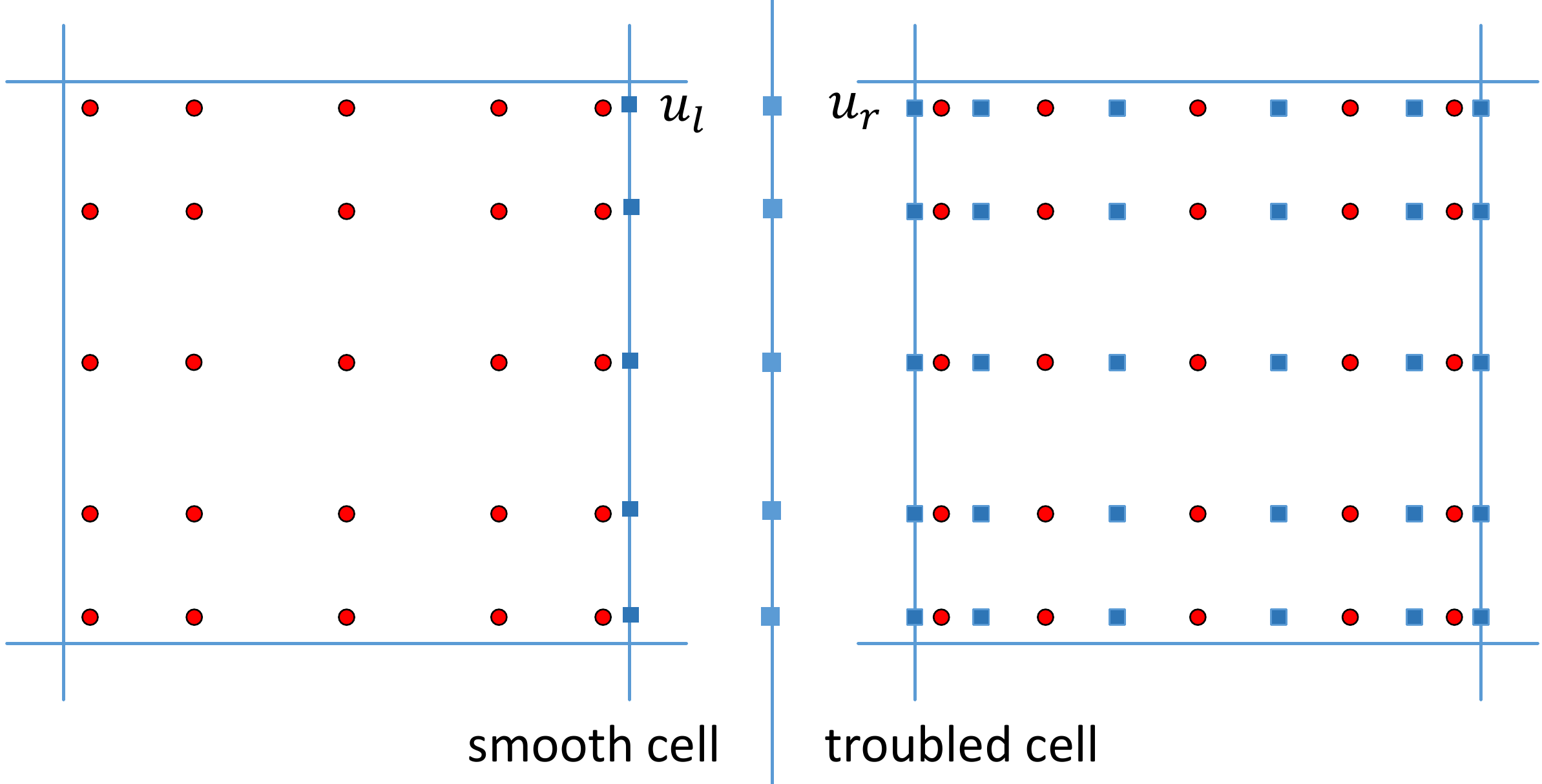}
			\caption{The cell interface between the smooth cell and the troubled cell. Flux points are marked in blue square and solution points are marked in red circle.}	
			\label{SchemeSwitch}	
		\end{figure}

		Like \cite{Sonntag2016}, we divide the CPR cell into subcells and apply CNNW2 on these subcells. The method can be interpreted as a blockwise CNNW2 method on unstructured quadrilateral blocks, where the blocks are coupled through Riemann solvers. As shown in Figure \ref{SchemeSwitch}, the left values of the Riemann flux on the cell interface are offered by the smooth cell, thus we compute the $u_l$ by Lagrange interpolation while the $u_r$ by nonuniform nonlinear weighted interpolation. Thus, the interpolation method to obtain the left and right value depending on the cell type, smooth cell or troubled cell. 
\section{Numerical investigation}\label{numerical results}
	\subsection{The isentropic vortex problem}
	In this subsection, we test the accuracy of CPR and CNNW2 scheme by solving 2D isentropic vortex problem \cite{Shu1997}. The initial condition is a mean flow with isotropic vortex perturbations. The mean flow is $\{\rho_{\infty},u_{\infty},v_{\infty},p_{\infty}\}=\{1,1,0,1\}$ , with $T_{\infty}=p_{\infty}/\rho_{\infty}$. 
	Consider the computational domain $\Omega=[-5,5]\times[-5,5]$, and the vortex is centered at $(x_c,y_c)=(0,0)$ with the following initial conditions
	\begin{equation}\begin{aligned}
	&\Delta u=-(y-y_c) \frac{\epsilon}{2 \pi} \exp \left(\frac{1-r^{2}}{2}\right), \quad \Delta v=(x-x_c) \frac{\epsilon}{2 \pi} \exp \left(\frac{1-r^{2}}{2}\right), \quad \\
	&\Delta T=-\frac{(\gamma-1) \epsilon}{8 \gamma \pi^{2}} \exp \left(1-r^{2}\right),
	\end{aligned}
	\end{equation}
	where $r=\sqrt{(x-x_c)^2+(y-y_c)^2}$ and the vortex strength $\epsilon=5$. With the state equation of ideal gas, $P=\rho R T,~P/\rho^{\gamma}=Constant$, the initial density and pressure follows the relation
	\begin{equation}
	P=P_{\infty}\left(\frac{T}{T_{\infty}}\right)^{\frac{\gamma}{\gamma-1}} \quad \rho=\rho_{\infty}\left(\frac{T}{T_{\infty}}\right)^{\frac{1}{\gamma-1}},
	\end{equation}
	so the initial condition is 
	\begin{equation}\label{state 1}
	(\rho, u, v,~p)=((T_{\infty}+\Delta T)^{\frac{\gamma}{\gamma-1}},~u_{\infty}+\Delta u,~v_{\infty}+\Delta v,~(T_{\infty}+\Delta T)^{\frac{1}{\gamma-1}}). 
	\end{equation}
	The problem is solved till $T$ = 0.2 with periodic boundary conditions.
	\begin{table}[H]\centering
		\caption{Accuracy test for 2D Euler equations}
		\label{CNNW2-order-t02}
		\begin{tabular}{cccccccc}
			\toprule
			&	Cells  & $L_1$ error &$L_1$ Order & $L_2$ error & $L_2$ Order& $L_{\infty}$ error& $L_{\infty}$ Order  \\ \midrule
	CNNW2 	&10$\times$10    & 2.86E-04   & --     & 7.16E-04   & --     & 5.72E-03     & --     \\
with no limiter	&20$\times$20    & 9.06E-05   & 1.66   & 2.32E-04   & 1.63   & 2.22E-03     & 1.37   \\
			&40$\times$40    & 2.58E-05   & 1.81   & 7.07E-05   & 1.71   & 8.34E-04     & 1.41   \\
			&80$\times$80    & 6.39E-06   & 2.01   & 1.72E-05   & 2.04   & 2.00E-04     & 2.06   \\ \hline
	CNNW2  	&10$\times$10    & 3.44E-04   & --     & 9.92E-04   & --     & 9.87E-03     & --     \\
with limiter&20$\times$20    & 1.01E-04   & 1.77   & 3.05E-04   & 1.70   & 4.73E-03     & 1.06   \\
			&40$\times$40    & 2.76E-05   & 1.87   & 8.39E-05   & 1.86   & 1.39E-03     & 1.77   \\
			&80$\times$80    & 6.89E-06   & 2.00   & 2.23E-05   & 1.91   & 4.40E-04     & 1.66   \\ \hline
	CPR (N=4) &10$\times$10    & 6.29E-06   & --     & 1.65E-05   & --     & 1.43E-04     & --      \\
			&20$\times$20    & 2.35E-07   & 4.75   & 6.49E-07   & 4.67   & 6.34E-06     & 4.50   \\
			&40$\times$40    & 1.08E-08   & 4.44   & 3.07E-08   & 4.40   & 3.46E-07     & 4.20   \\
			&80$\times$80    & 5.82E-10   & 4.22   & 1.58E-09   & 4.28   & 1.62E-08     & 4.41   \\ \hline
		\end{tabular}
	\end{table}
	$L_{1},~L_{2}$ and $L_{\infty}$ error norms for the density at the final time are computed by 
	\begin{equation}
	\text {error}_{L_{1}}=\frac{\sum_{i=1}^{N_p}|U_{i}^{h}-U_{i}|}{N_p},~
	\text {error}_{L_{2}}=\sqrt{\frac{\sum_{i=1}^{N_p}\left(U_{i}^{h}-U_{i}\right)^{2}}{N_p}},~
	\text {error}_{L_{\infty}}=max\{U_{i}^{h}-U_{i}\}_{i=1}^{N_p},	
	\end{equation}
	where $N_p$ is the number of the solution points, 
	$U_{i}^{h}$ and $U_{i}$ are the numerical solution and 
	the exact solution at the i-th solution point,
	and $N_p$ is the number of the solution points. 
	The errors and the rates of convergence are shown in Table \ref{CNNW2-order-t02}.
	The CNNW2 scheme with no limiter gets a convergence rate of about 2, which coincides with the designed order of accuracy. And the accuracy of the CNNW2 scheme with limiter is lower than that with no limiter because the limiter plays a role at the extreme points and reduces the accuracy. 
	We study the subcell limiting strategy of the high-order CPR method which achieves fifth-order in our procedure. However, the numerical test for nonlinear conservation laws reveals that there is a slight accuracy loss (0.5 $\sim$ 1 order) with LP (Lagrange polynomial) approach to compute the divergence of the flux vector \cite{Wang2009}. 
	\subsection{Sod and Lax shock tube problems}
\begin{figure}[htbp]\centering
		\includegraphics[scale=0.7]{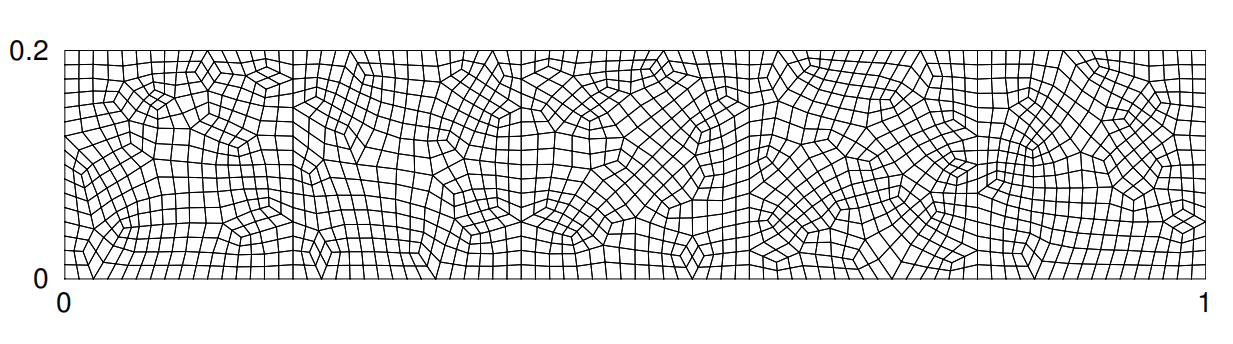}
		\caption{The computational meshes of shock tube problem in total 1604 cells.}\label{sod lax mesh}		
	\end{figure}

	The planar Sod and the classical Lax shock tube problems are solved by CPR-CNNW2 on 2D unstructured quadrilateral meshes, as shown in Figure \ref{sod lax mesh}. These problems, have analytical solutions \cite{Toro2009}, are used to assess the ability of the proposed numerical method in capturing one-dimensional discontinuous flows. Dirichlet boundary conditions and periodic boundary conditions are imposed in the horizontal and vertical directions, respectively. 
	\subsubsection{Sod problem}
	The classical Sod problem is initialized by a discontinuity located in the middle of the computational domain $\Omega=[0, 1]$ \cite{Sod1978}. The initial condition of Sod problem is
	\begin{equation}
	(\rho,~u,~v,~p)=\left\{\begin{array}{ll}
	(0.125,~ 0,~ 0,~ 0.1), & 0 \leq x<0.5, \\
	(1.0,~0,~0,~ 1.0), & 0.5 \leq x \leq 1.	\end{array}\right.
	\end{equation}
	The problem is solved till T = 0.2 by hybrid CPR-CNNW2.\par
\begin{figure}[H]
		\centering
		\subfigbottomskip=2pt
		\subfigcapskip=-5pt
		\subfigure[]{\label{sod rho}
			\includegraphics[width=0.45\linewidth]{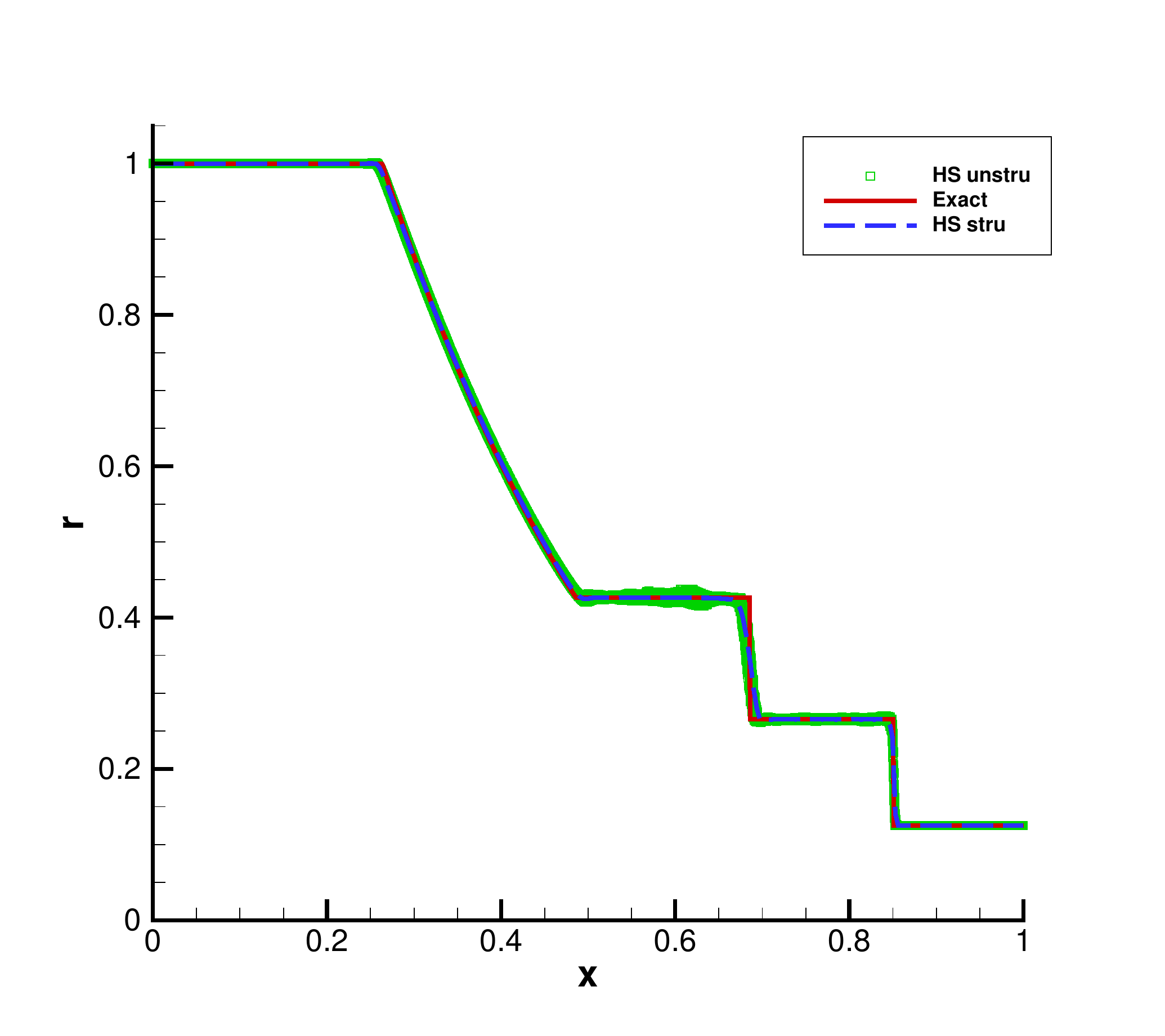}}
		\subfigure[]{\label{sod beta}
			\includegraphics[width=0.45\linewidth]{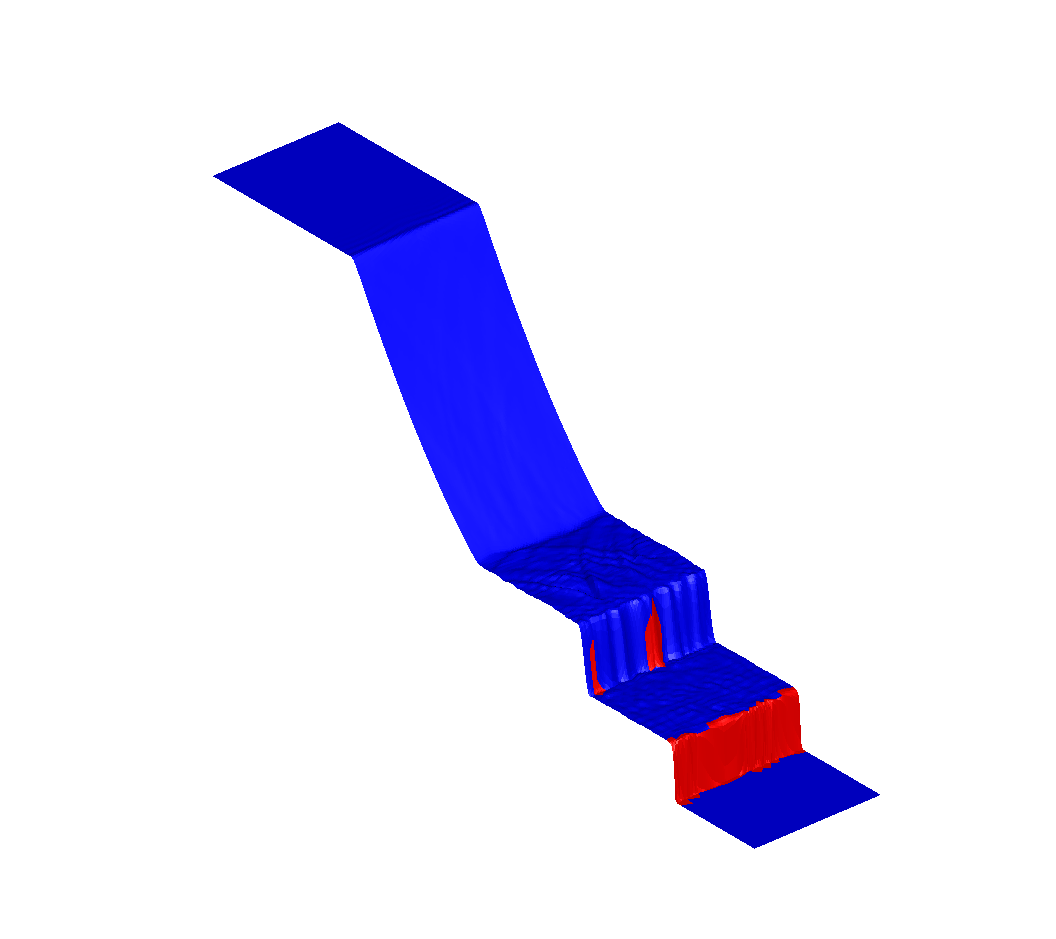}}
		\\
		\subfigure[]{\label{sod rho planar contour map}
			\includegraphics[width=0.45\linewidth]{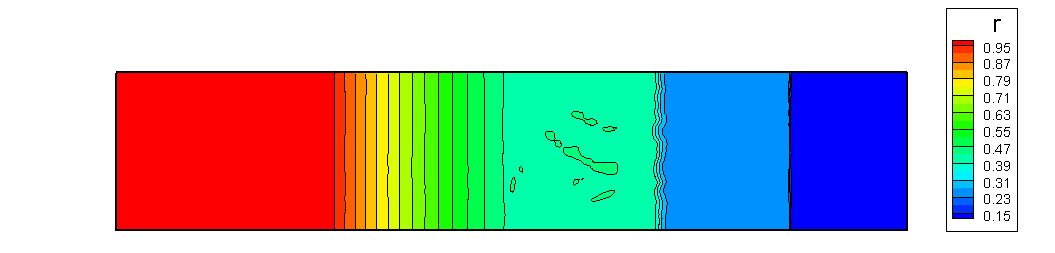}}
		\subfigure[]{\label{sod beta planar}
			\includegraphics[width=0.45\linewidth]{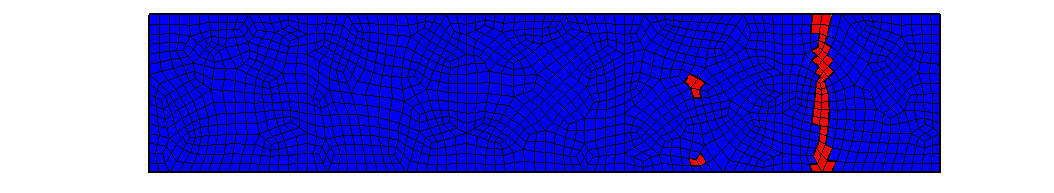}}
		\caption{Results calculated by the hybrid CPR-CNNW2 scheme at T = 0.2. (a) Density $\rho$. The result is the projection on the x-z plane. The result computed on unstructured meshes mark in green square, and on structured meshes in blue dash. (b)The distribution of troubled cells: the red cells. (c) Density from 0.15 to 0.95 with 21 contours. (d) The planar distribution of troubled cells.}
	\end{figure}

	Figure \ref{sod rho} shows that the subcell limiting strategy based on CNNW2 can capture shock waves effectively. The troubled cells accounts for 2.56 \% of the total number of cells and mainly concentrated near the contact discontinuity (x =0.7) and the shock (x = 0.85), as shown in Figure \ref{sod beta} and \ref{sod beta planar}. 
	Compared with the result on structured meshes( $90\times18$), the result on unstructured meshes have slight fluctuations as shown in Figure \ref{sod rho}, which is due to the irregularity of the unstructured meshes \cite{Clain2011}. 
	\subsubsection{Lax problem}
	Consider the Lax problem with initial condition 
	\begin{equation}
	(\rho,~u,~v,~p)=\left\{\begin{array}{ll}
	(0.445,~ 0.698,~ 0,~ 3.528), & 0 \leq x<0.5, \\
	(0.5,~0,~0,~ 0.571), & 0.5 \leq x \leq 1.	\end{array}\right.
	\end{equation}
	This problem is solved till T = 0.14.
	\begin{figure}[H]
		\centering
		\subfigbottomskip=2pt
		\subfigcapskip=-5pt
		\subfigure[]{\label{lax rho}
			\includegraphics[width=0.45\linewidth]{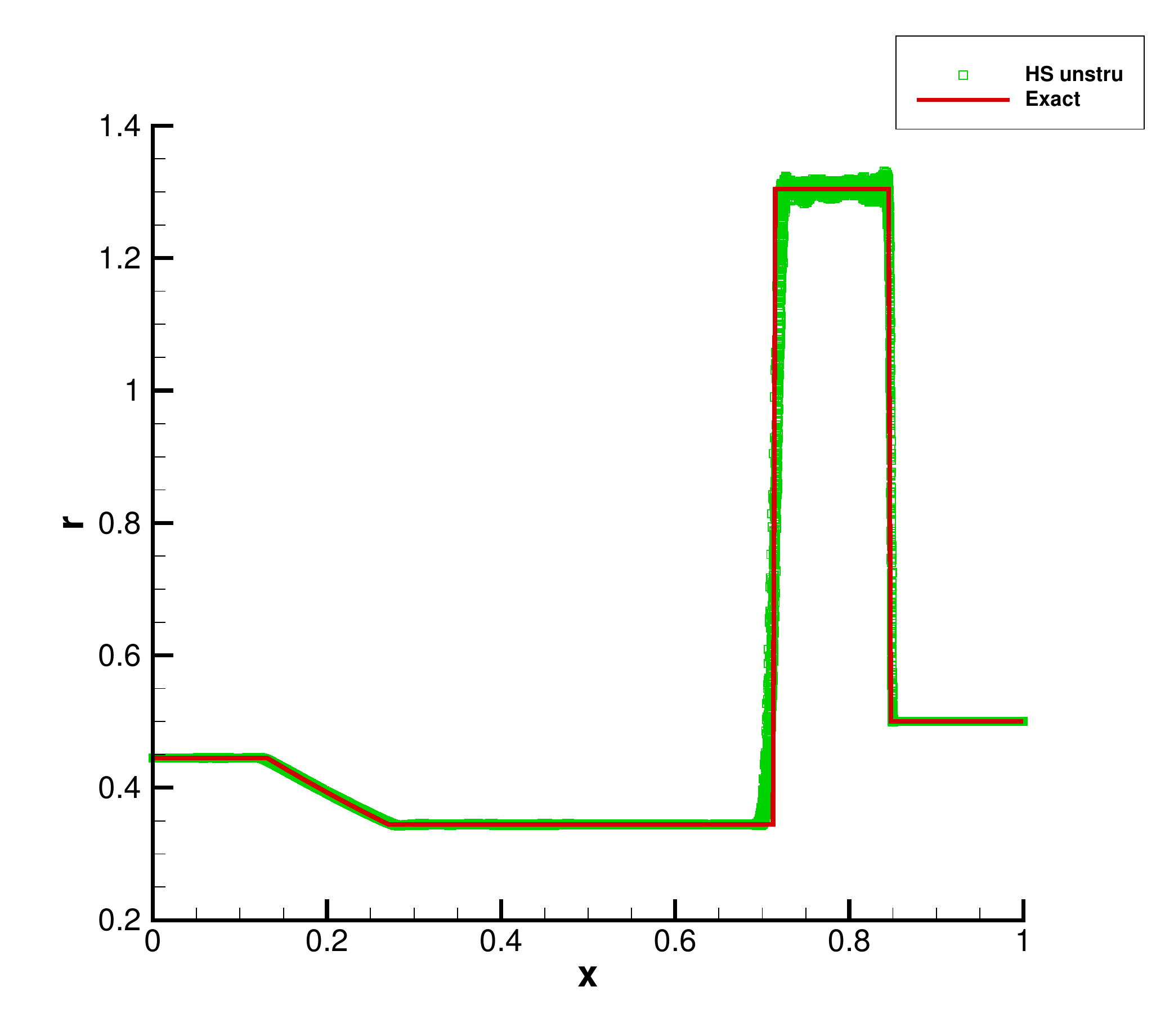}}
		\subfigure[]{\label{lax beta}
			\includegraphics[width=0.45\linewidth]{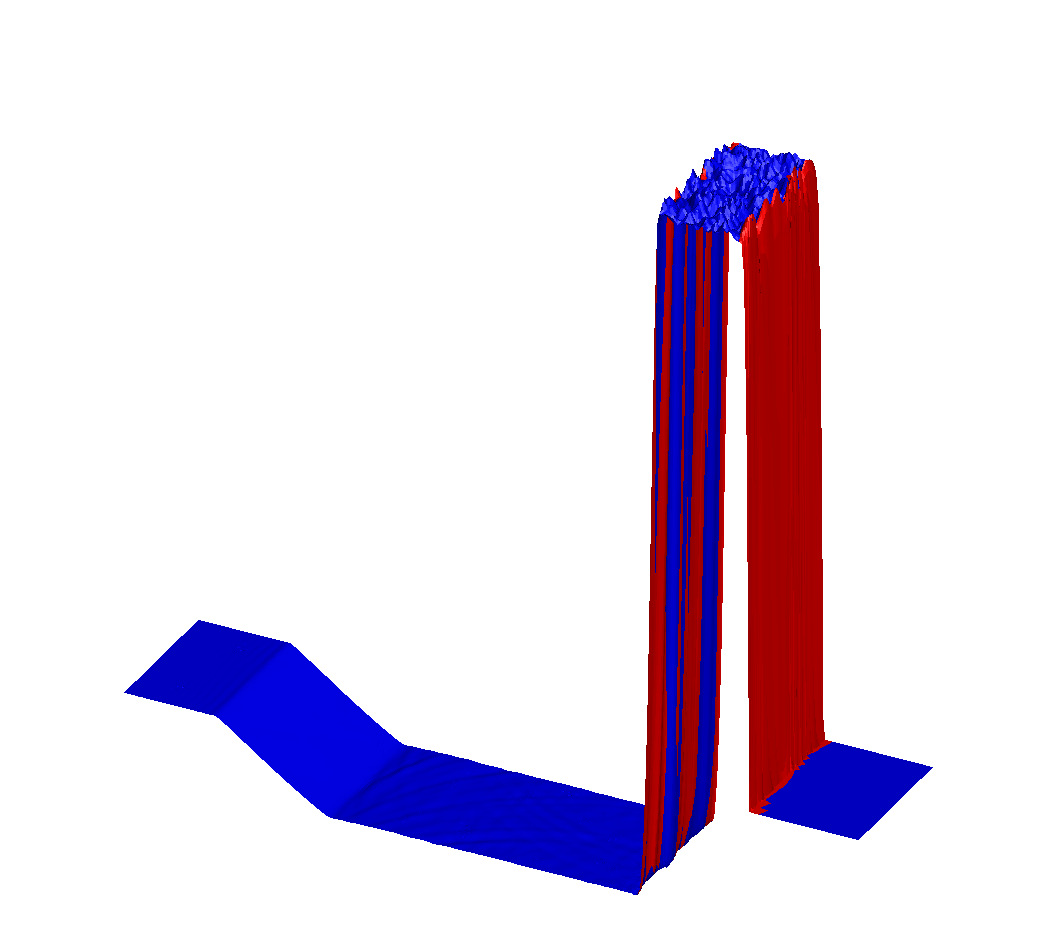}}
		\\
		\subfigure[]{\label{lax rho planar contour map}
			\includegraphics[width=0.45\linewidth]{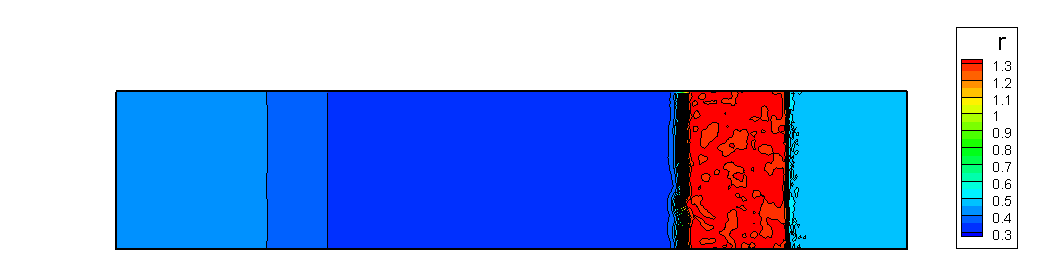}}
		\subfigure[]{\label{lax beta planar}
			\includegraphics[width=0.45\linewidth]{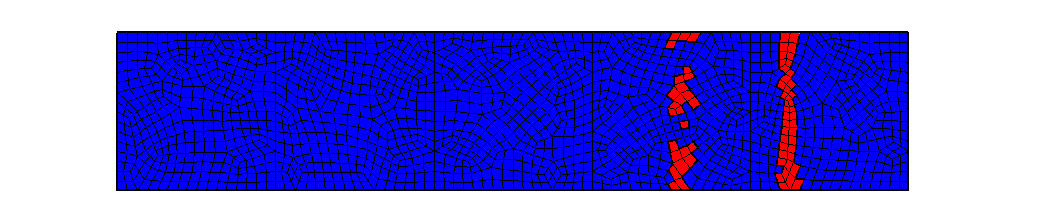}}
		\caption{Results calculated by the hybrid CPR-CNNW2 scheme. T = 0.14. (a) Density $\rho$. The result is the projection on the x-z plane. (b) The distribution of troubled cells: the red cells. (c) Density from 0.3 to 1.3 with 21 contours. (d) The planar distribution of troubled cells.}
	\end{figure}
	From Figure \ref{lax rho}, we can see that at x = 0.72 and x = 0.85, shocks and discontinuities are all captured and the troubled cells take a percentage of 4.11 \% of the total cells.
	The troubled cells marked red are shown in Figure \ref{lax beta} and \ref{lax beta planar}.
	\subsection{Shu-Osher problem}
	This problem was proposed by Shu and Osher in \cite{Shu1989} to study the ability in capturing discontinuities and high-frequency waves interacting with the discontinuities. The initial condition is divided into a high-pressure part on the left side and a sinusoidal density wave on the right side. We solve the Shu-Osher problem with the initial condition
	\begin{equation}
	(\rho,~ u,~ v,~ p)=\left\{\begin{array}{lc}
	(3.857143,~2.629369,~0,~10.33333), & 0 \leq x<0.1, \\
	(1.0+0.2 \sin (50 x),~ 0,~0,~1.0), & 0.1<x \leq 1.0.
	\end{array}\right.
	\end{equation}
	till T = 0.18 on the unstructured meshes. The boundary conditions and the computational meshes are the same as the Sod problem. 
	\begin{figure}[H]
		\centering
		\subfigbottomskip=2pt
		\subfigcapskip=-5pt
		\subfigure[]{\label{shu-osher rho}
			\includegraphics[width=0.45\linewidth]{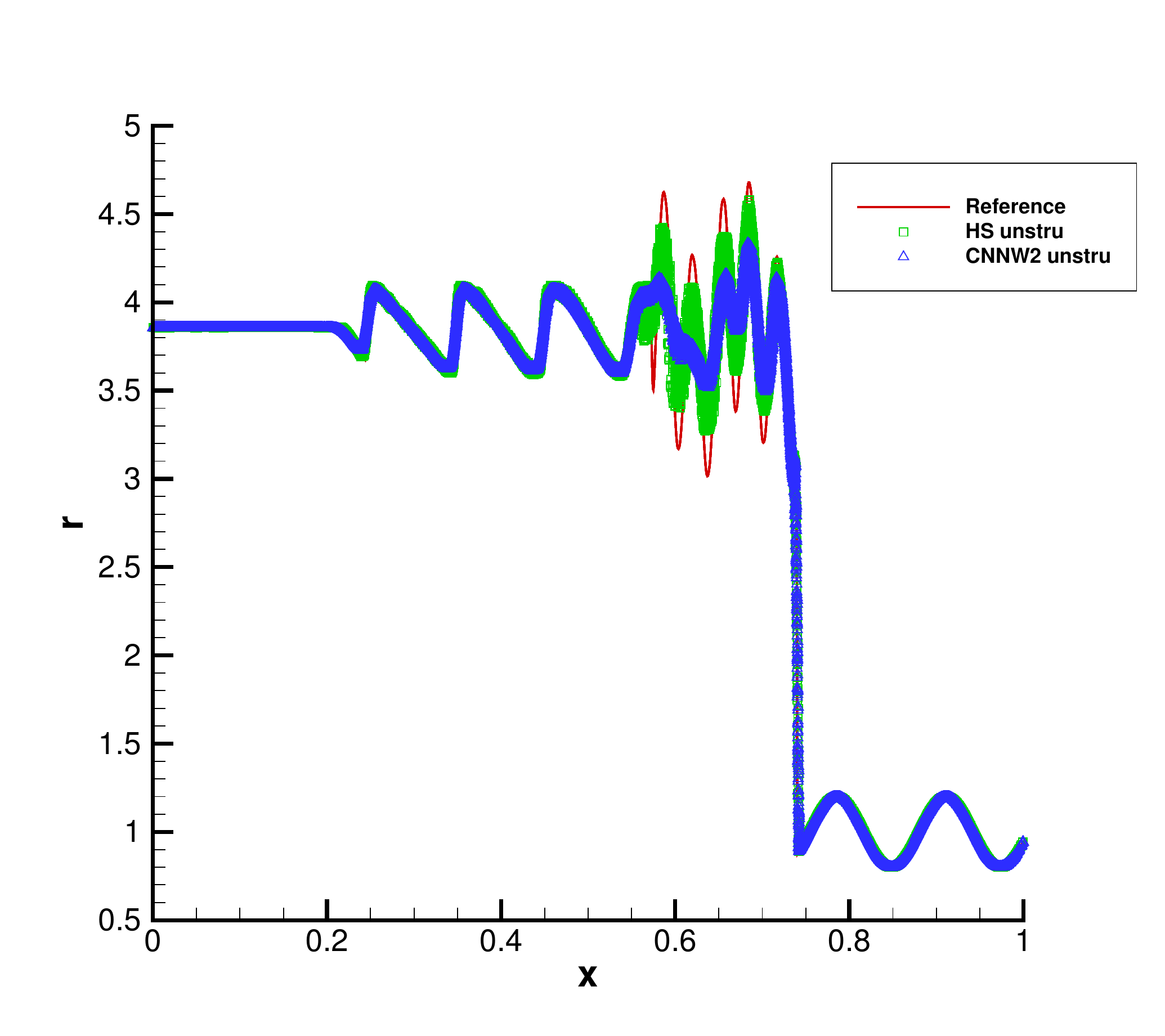}}
		\subfigure[]{\label{shu-osher beta}
			\includegraphics[width=0.45\linewidth]{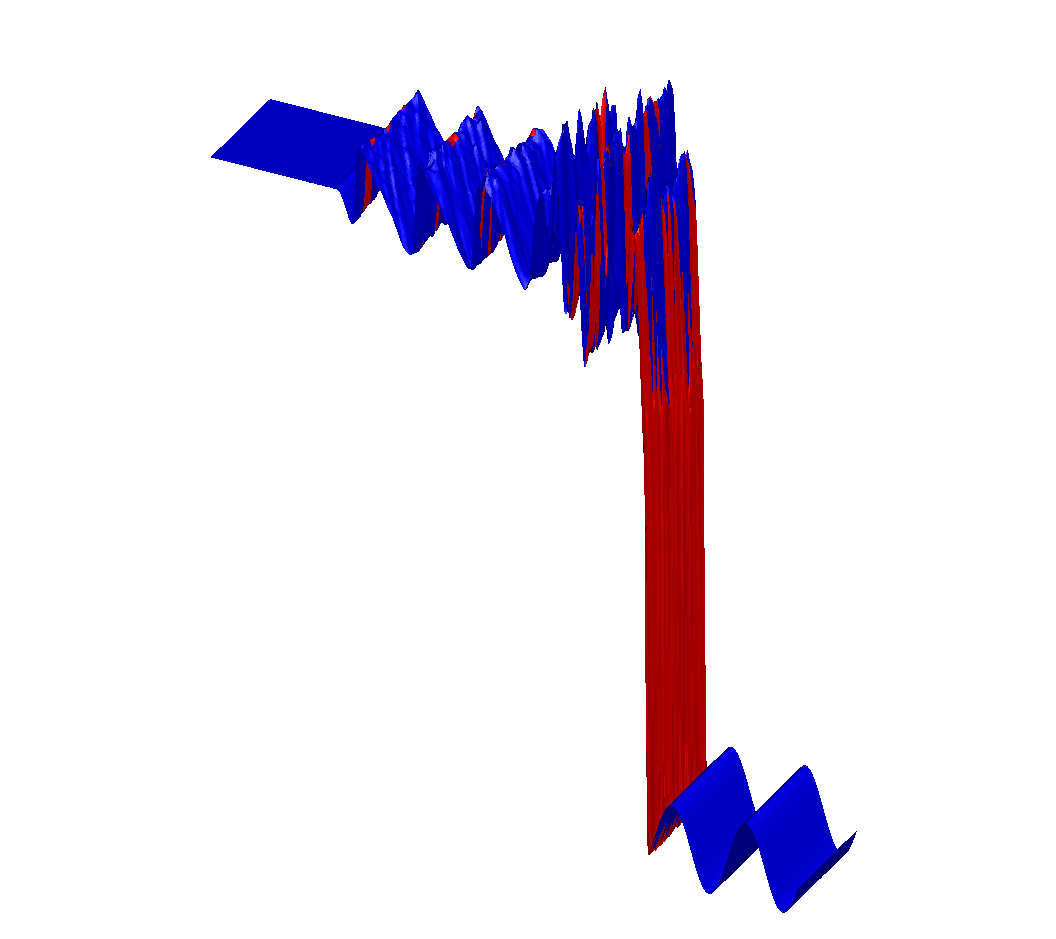}}
		\\
		\subfigure[]{\label{shu-osher rho planar}
			\includegraphics[width=0.45\linewidth]{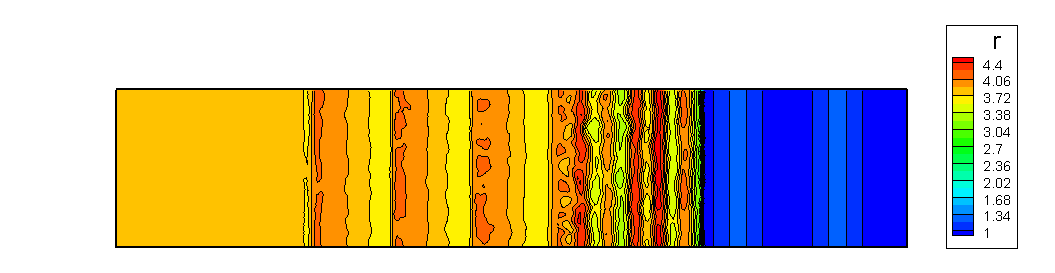}}
		\subfigure[]{\label{shu-osher beta planar}
			\includegraphics[width=0.45\linewidth]{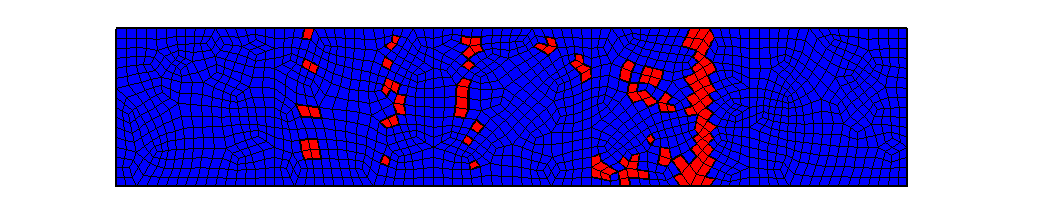}}
		\caption{Results calculated by the hybrid CPR-CNNW2 scheme at T = 0.18. (a) Density $\rho$. The result is the projection on the x-z plane. The result computed on unstructured meshes with the hybrid scheme is marked in green square, with CNNW2 scheme in blue delta. (b) The distribution of troubled cells: the red cells. (c) Density from 1.0 to 4.4 with 21 contours. (d) The planar distribution of troubled cells.}
		\label{shu-osher rho beta}
	\end{figure}
	\begin{figure}[htbp]\centering
		\includegraphics[scale=0.4]{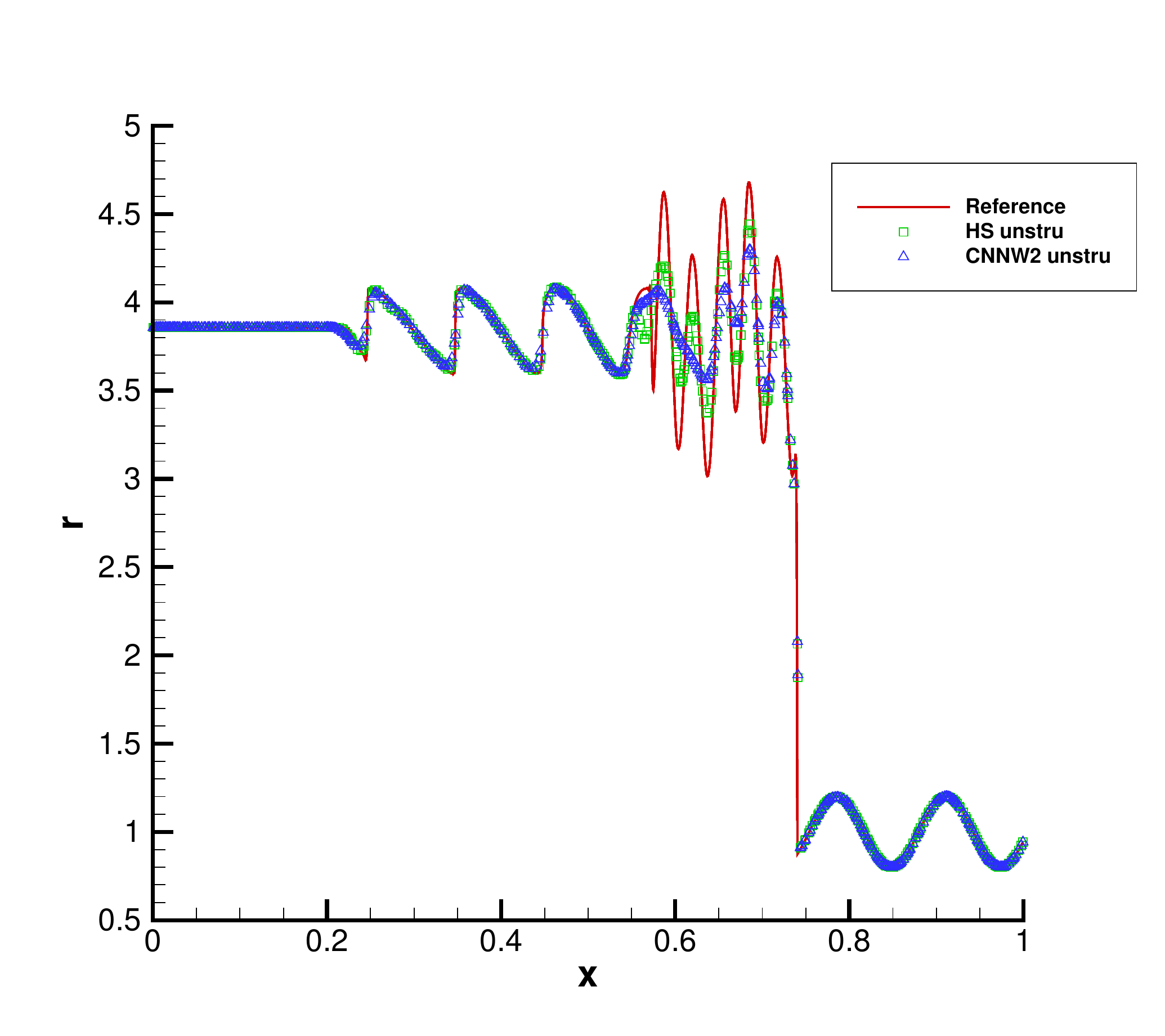}
		\caption{The density slice at y = 0.1.}\label{shu-osher-slice}		
	\end{figure}
	As shown in Figure \ref{shu-osher rho beta} and Figure \ref{shu-osher-slice}, the reference solution is obtained by the hybrid CPR-CNNW2 scheme on structured meshes (400 $\times$ 2), with DOFs = 2000 in the x-direction. And the hybrid scheme can capture shocks and discontinuities well on unstructured meshes while 6.61 \% of the total cells are troubled cells. Moreover, comparisons between the hybrid scheme and CNNW2 are made, and the results illustrate that the former has a higher resolution than the latter one.
	\subsection{2D Riemann problem}
	The two-dimensional Riemann problem is divided into 19 typical types by the structure of the solution \cite{Lax1998}. We test the performance of the subcell limiting strategy on the shock wave interaction problem.
	With
	\begin{equation}
	(\rho,~ u,~ v,~ p)=\left\{\begin{array}{ll}
	(1.5,~0,~0,~1.5), & \text { if } 0.8<x<1.0 \text { and } 0.8<y<1.0, \\
	(0.5323,~1.206,~0,~0.3), & \text { if } 0<x<0.8 \text { and } 0.8<y<1.0, \\
	(0.138,~1.206,~1.206,~0.029), & \text { if } 0<x<0.8 \text { and } 0<y<0.8, \\
	(0.5323,~0,~1.206,~0.3), & \text { if } 0.8<x<1.0 \text { and } 0<y<0.8,
	\end{array}\right.
	\end{equation}
	and Dirichlet boundary conditions, this problem is run until T=0.8 on the computational domain $[0, 1]\times[0, 1]$.	
\par 
	\begin{figure}[htbp]\centering
		\subfigure[]{\label{2D Riemann rho}
			\includegraphics[width=0.48\linewidth]{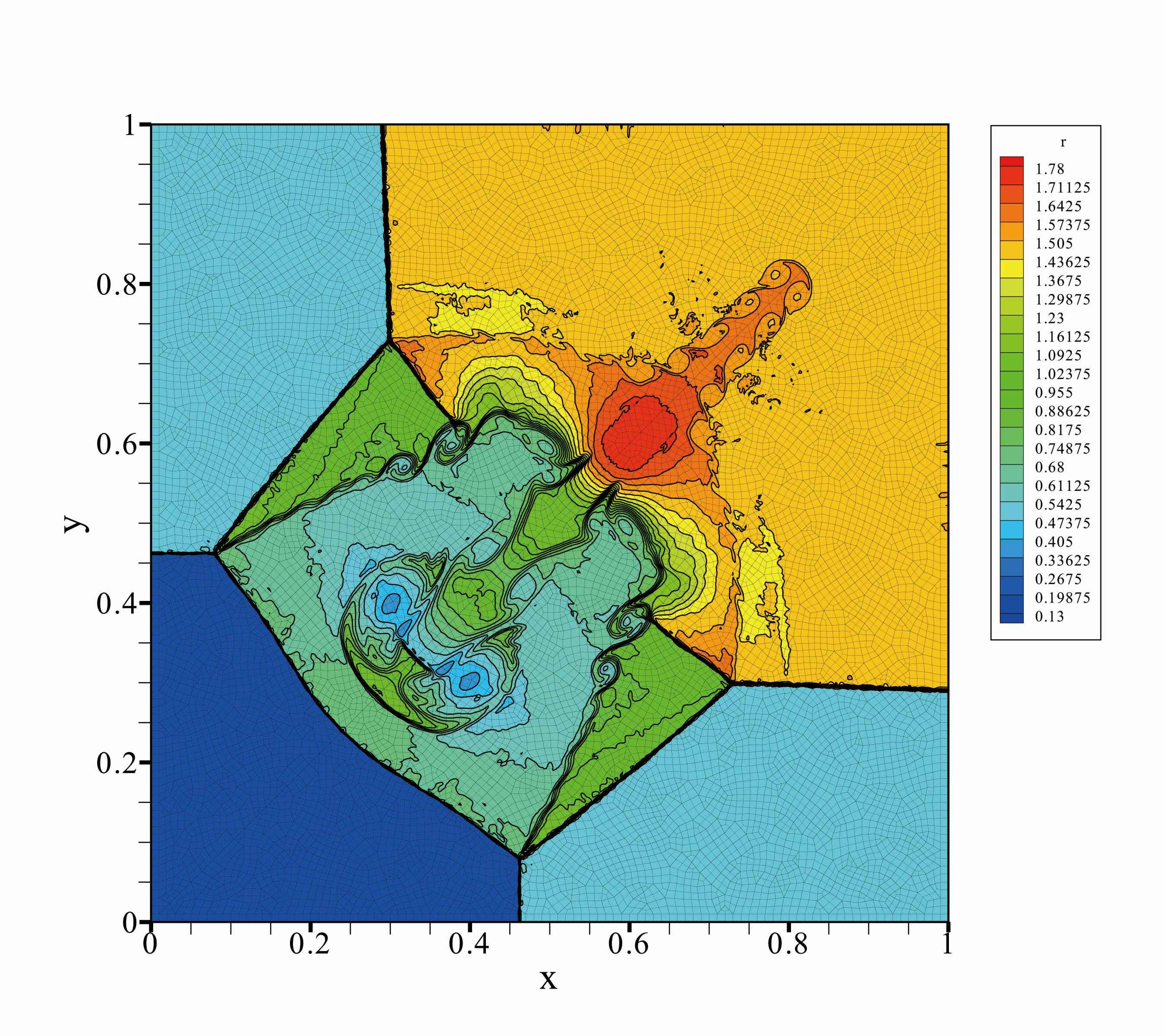}}
		\subfigure[]{\label{2D Riemann beta}
			\includegraphics[width=0.48\linewidth]{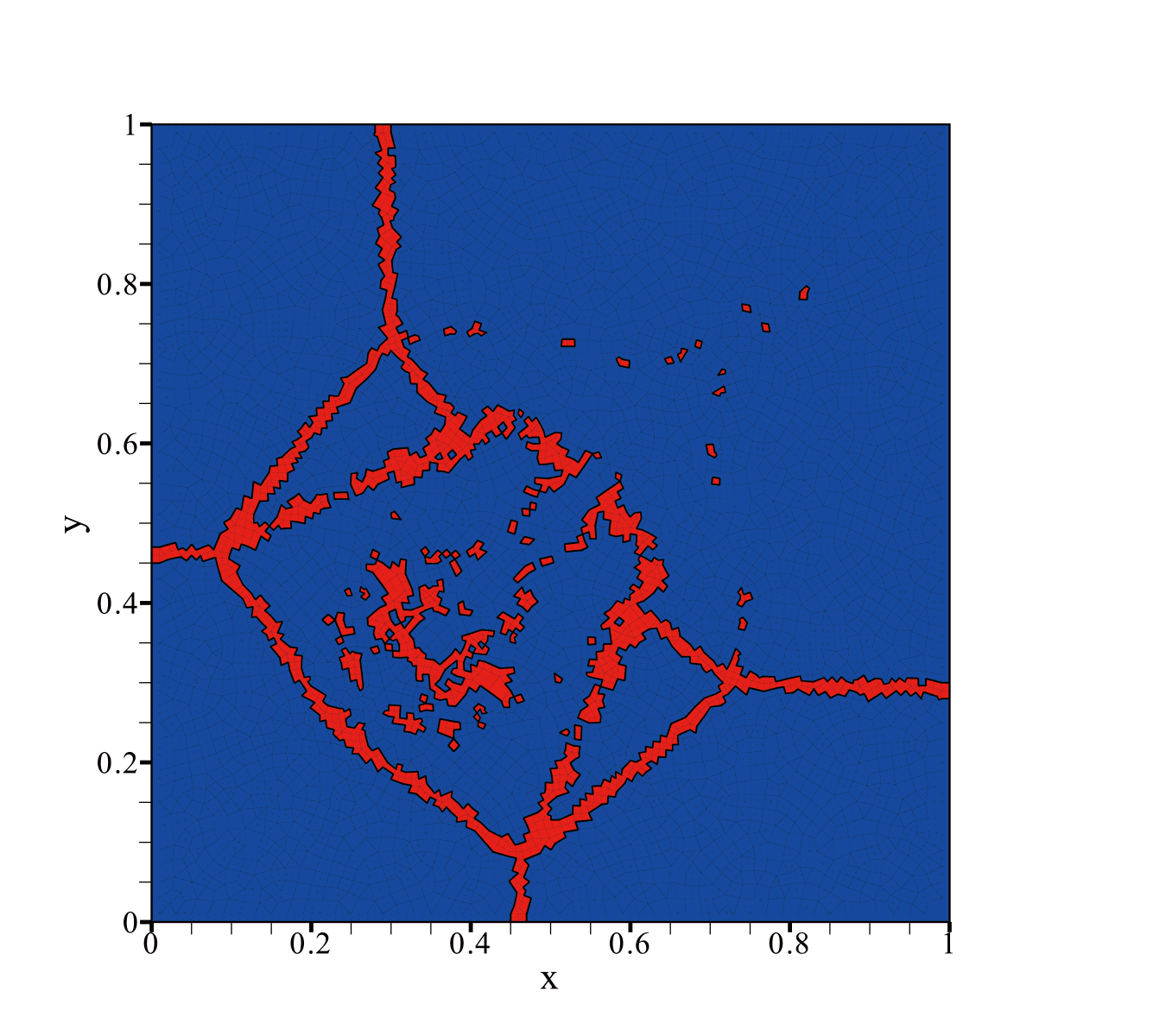}}
		\caption{Results calculated by the hybrid scheme at T = 0.8 with 14494 cells. (a) Density $\rho$ from 0.13 to 1.78 with 25 contours. (b) Distribution of troubled cells: the red cells.}	
	\end{figure}
	The results of the hybrid scheme in Figure \ref{2D Riemann rho} show that small-scale structures are well captured. As shown in Figure \ref{2D Riemann beta}, we can see that only 8.62 \% of the computational domain is calculated by CNNW2 and most cells are captured by CPR, which contributes to the high efficiency of the hybrid scheme. 
	\subsection{Double Mach reflection problem}
	2D double Mach reflection problem with a strong shock is proposed by Woodward and Colella \cite{Paul1984} to test the robustness of the high-resolution schemes. The computational domain is $[0, 4]\times[0, 1]$. This test problem involves a Mach 10 shock that is initially set up at x = 1/6 on the lower boundary and shapes a $60 \degree$ ramp with the x-axis. Specifying inflow boundary condition at the bottom boundary for x = [0, 1/6] and slip wall boundary condition for x = [1/6, 4]. For the upper boundary (y = 1), a time-dependent boundary based on the analytical propagation speed of the oblique shock is imposed. Inflow boundary conditions and outflow boundary conditions are set at the left and right boundary respectively. The initial condition is
	\begin{equation}
	(\rho,~ u,~ v,~ p)=\left\{\begin{array}{ll}
	(1.4,~0.0,~0.0,~1.0), & \text { if } y<\sqrt{3}\left(x-\frac{1}{6}\right), \\
	(8.0,~7.145,~-4.125,~116.5), & \text { if } y \geq \sqrt{3}\left(x-\frac{1}{6}\right) .
	\end{array}\right.
	\end{equation}
	and the final computing time is set at T = 0.2.
	\begin{figure}[H]
		\centering
		\subfigbottomskip=2pt
		\subfigcapskip=-5pt
		\subfigure[]{\label{doubleMach rho}
			\includegraphics[width=0.8\linewidth]{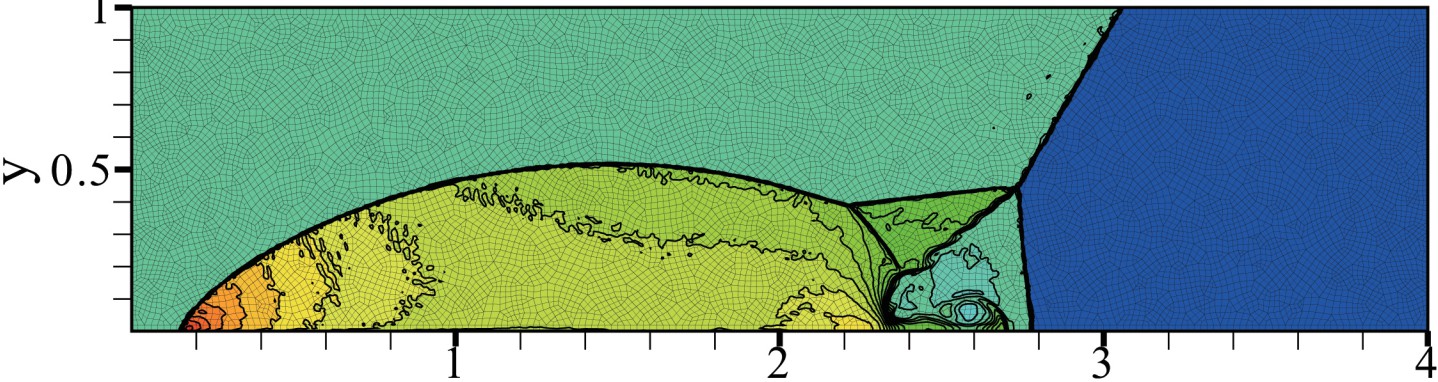}}
		\\
		\subfigure[]{\label{doubleMach beta}
			\includegraphics[width=0.8\linewidth]{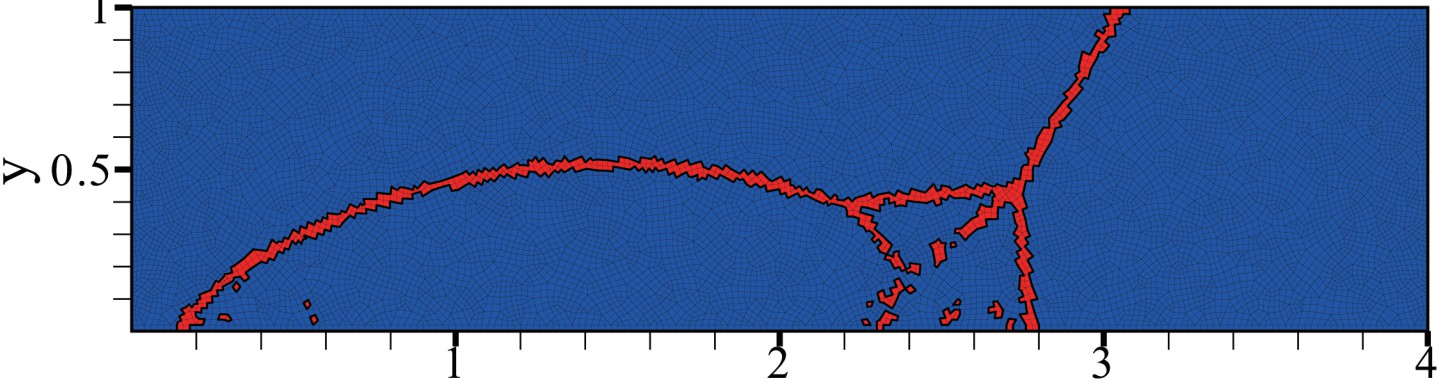}}
		\caption{Results calculated by the hybrid CPR-CNNW2 scheme at T = 0.2 with 14464 cells. (a) Density $\rho$ from 1.5 to 21.7 with 30 contours. (b) Distribution of troubled cells: the red cells.}
	\end{figure}
	As shown in Figure \ref{doubleMach rho}, the hybrid CPR-CNNW2 scheme works well for this test case, and essential small- and large-scale structures of the flow field are captured. There are only 3.84 \% of cells that need to use the subcell limiting, which is marked in red in Figure \ref{doubleMach beta}. 
	\subsection{The shock-isentropic vortex interaction problem}\label{The ShockVortex Interaction Problem}
\begin{figure}[htbp]\centering
		\includegraphics[scale=0.45]{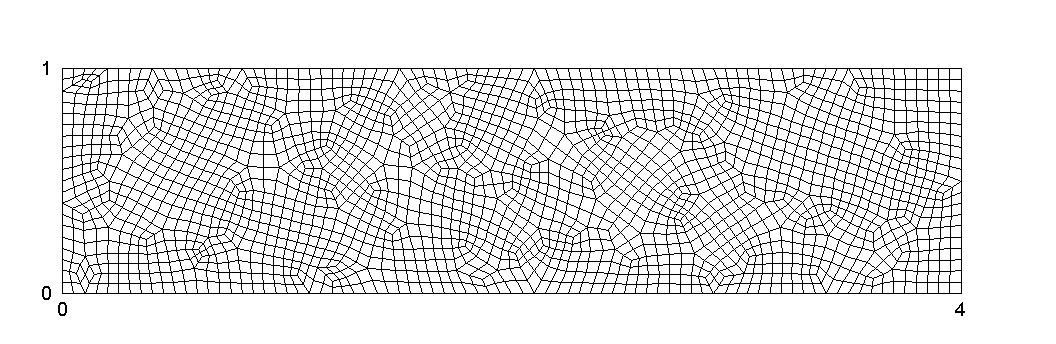}
		\caption{The computational meshes of the shock-vortex problem with 2004 cells.}	\label{shock vortex mesh}		
	\end{figure}

	The next problem is the shock-vortex interaction, which is used by Jiang and Shu \cite{Jiang1996} as a test problem. Here we consider the computational domain $\Omega = [0,4]\times[0,1]$, which shown as Figure \ref{shock vortex mesh}. A stationary Mach 1.1 shock is located at x = 2.0 and an isentropic vortex is superposed to the flow left to the shock and centers at $(x_c, y_c) = (0.25,0.5)$. The vortex perturbation is given as follows
	\begin{equation}
	\begin{aligned}
	&(\Delta u, \Delta v)= \epsilon \tau \exp (\alpha(1-{\tau}^{2}))(sin\theta,-cos\theta), \\
	&\Delta T=-\frac{(\gamma-1) \epsilon^2}{4 \alpha \gamma } \exp(2\alpha(1-\tau^{2})),
	\end{aligned}
	\end{equation}
	where $\tau=r/r_c,~r=\sqrt{(x-x_c)^2+(y-y_c)^2}$, $r_c =0.05$ is the critical radius of the vortex,  $\epsilon=0.3$ is the vortex strength, and $\alpha=0.204$ is the controlling factor of the vortex size. 	
	The flows of the left ($x<2.0$) and right ($x>2.0$) of the shock, are uniform,  are denoted by 1 and 2, respectively. 
	With 
	\begin{equation}
	(\bar{\rho}_1,~\bar{u}_1,~\bar{v}_1,~\bar{p}_1)=(1.0,~\sqrt{\gamma}\bar{M}_1,~0.0,~1.0),
	\end{equation}
	the initial values of region 2 can be computed through the Rankine-Hugoniot conditions \cite{leveque2002}. 
	As the perturbation is introduced into the flow to the left of the shock, then the initial values of region 1 $(\rho_1,u_1,v_1,p_1)$ can be solved with the state equation of ideal gas \cite{Jiang1996}. 
\par
	The left and right boundaries are implemented as zero-order extrapolation boundary conditions and the top and bottom boundaries are implemented as slip-wall boundary conditions \cite{Jiang1996}.
	\begin{figure}[htbp]\centering
		\includegraphics[scale=0.5]{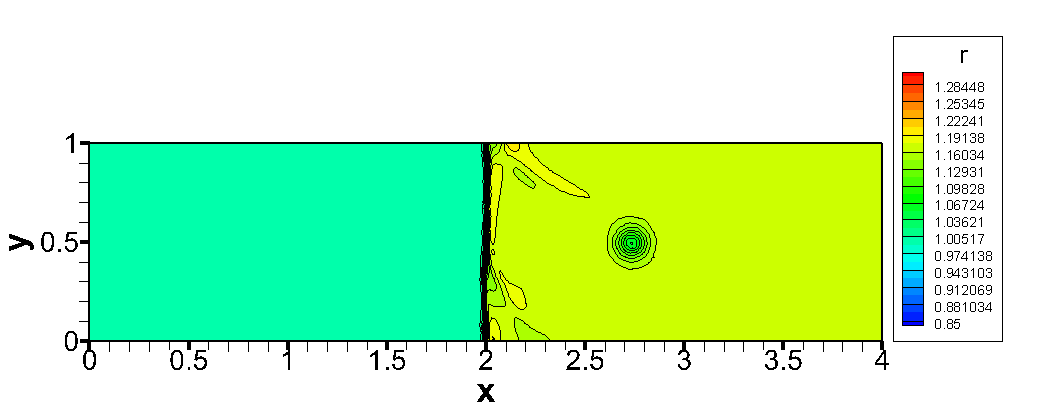}
		\caption{Results calculated by CNNW2 scheme at T = 2.0. Density $\rho$ from 0.85 to 1.3 with 30 contours.}	\label{CNNW2 shockvortex}		
	\end{figure}
	\begin{figure}[H]
		\centering
		\subfigbottomskip=2pt 
		\subfigcapskip=-5pt
		\subfigure[]{\label{shockvortex-hybrid}
			\includegraphics[width=0.95\linewidth]{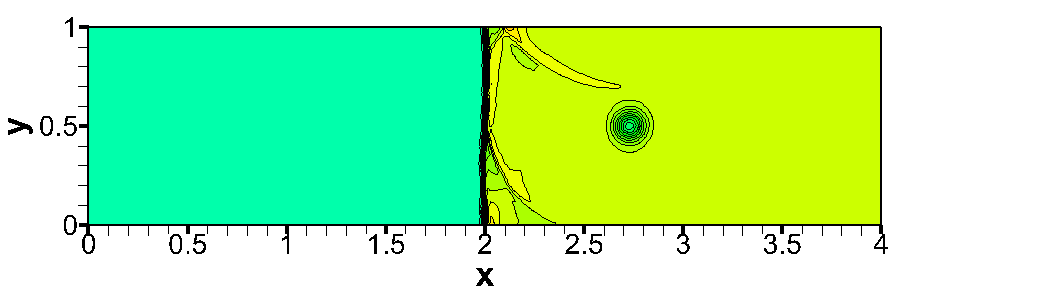}}
		\\
		\subfigure[]{\label{shockvortex-hybrid beta}
			\includegraphics[width=0.95\linewidth]{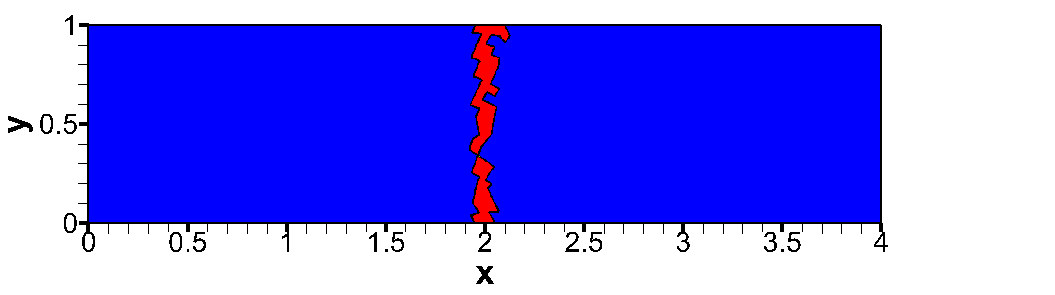}}
		\caption{Results calculated by the hybrid CPR-CNNW2 scheme at T = 2.0. (a) Density $\rho$ from 0.85 to 1.3 with 30 contours. (b) The distribution of troubled cells: the red area represents troubled cells, and the blue area represents smooth cells.}
	\end{figure}
	\begin{figure}[htbp]\centering
		\includegraphics[scale=0.6]{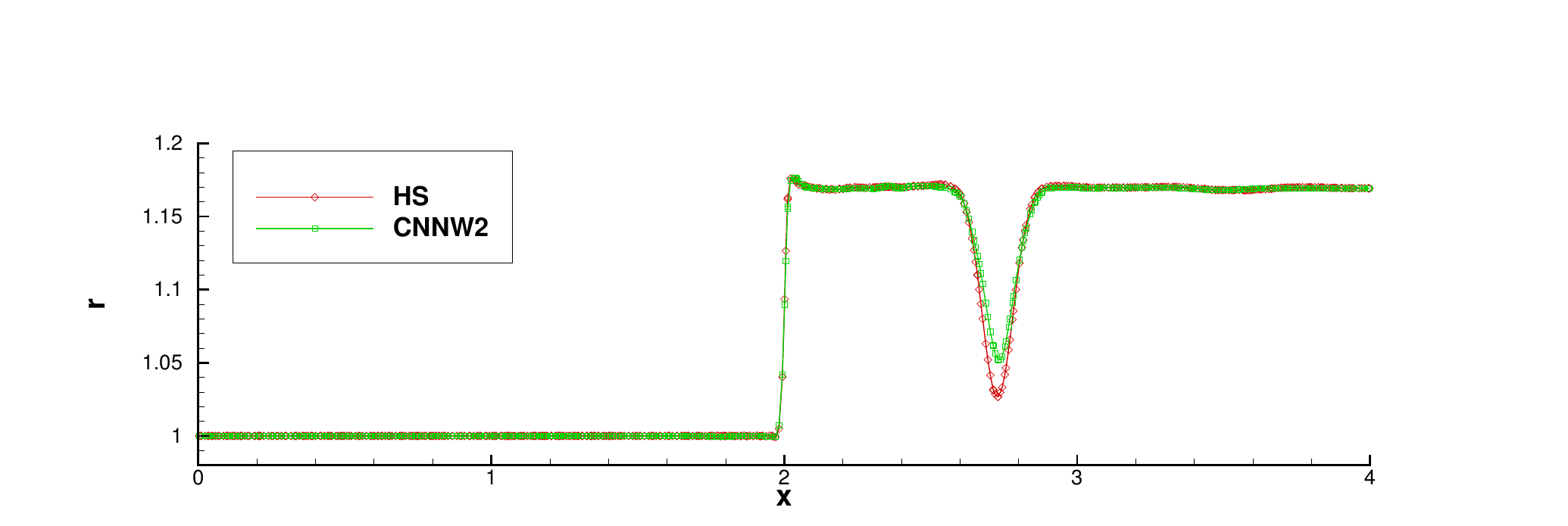}
		\caption{Results calculated by the hybrid CPR-CNNW2 scheme and CNNW2 scheme at t = 2.0. Slice at y = 0.5. }	\label{shockvortex-rho-compare}		
	\end{figure}
	The computation results are consistent with those in \cite{Jiang1996}. As shown in Figure \ref{CNNW2 shockvortex} and Figure \ref{shockvortex-hybrid}, the result of the hybrid scheme has a higher resolution than CNNW2, because most regions are computed by the CPR method, as the blue area in Figure \ref{shockvortex-hybrid beta}.  And from the slice at y = 0.5 shown in Figure \ref{shockvortex-rho-compare}, it is clearly illustrated that the CNNW2 scheme has more dissipation.
	\begin{figure}[H]
		\centering 
		\subfigbottomskip=2pt
		\subfigcapskip=-5pt
		\subfigure[]{\label{strongshockvortex-rho}
			\includegraphics[width=0.7\linewidth]{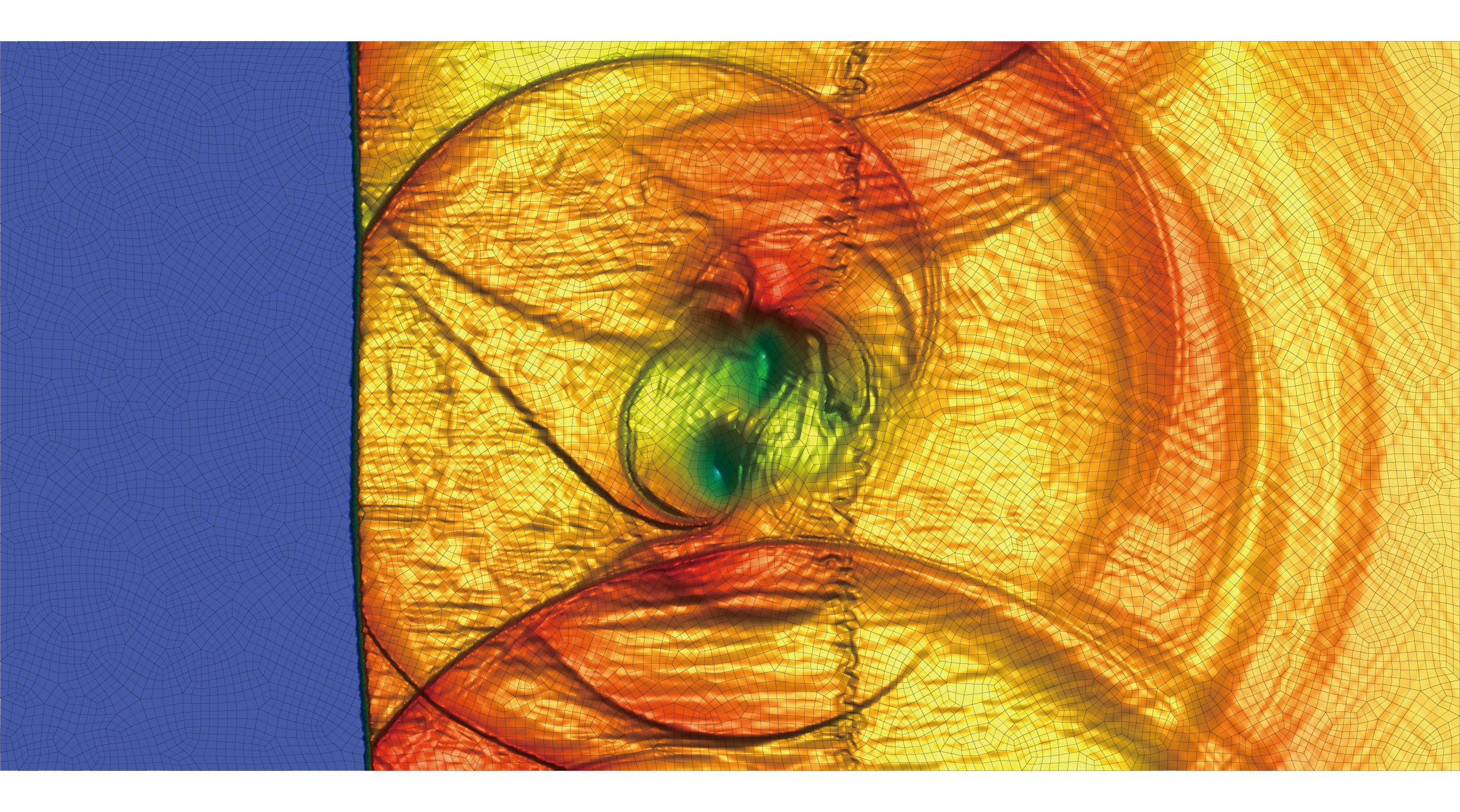}}
		\\
		\subfigure[]{\label{strongshockvortex-beta}
			\includegraphics[width=0.7\linewidth]{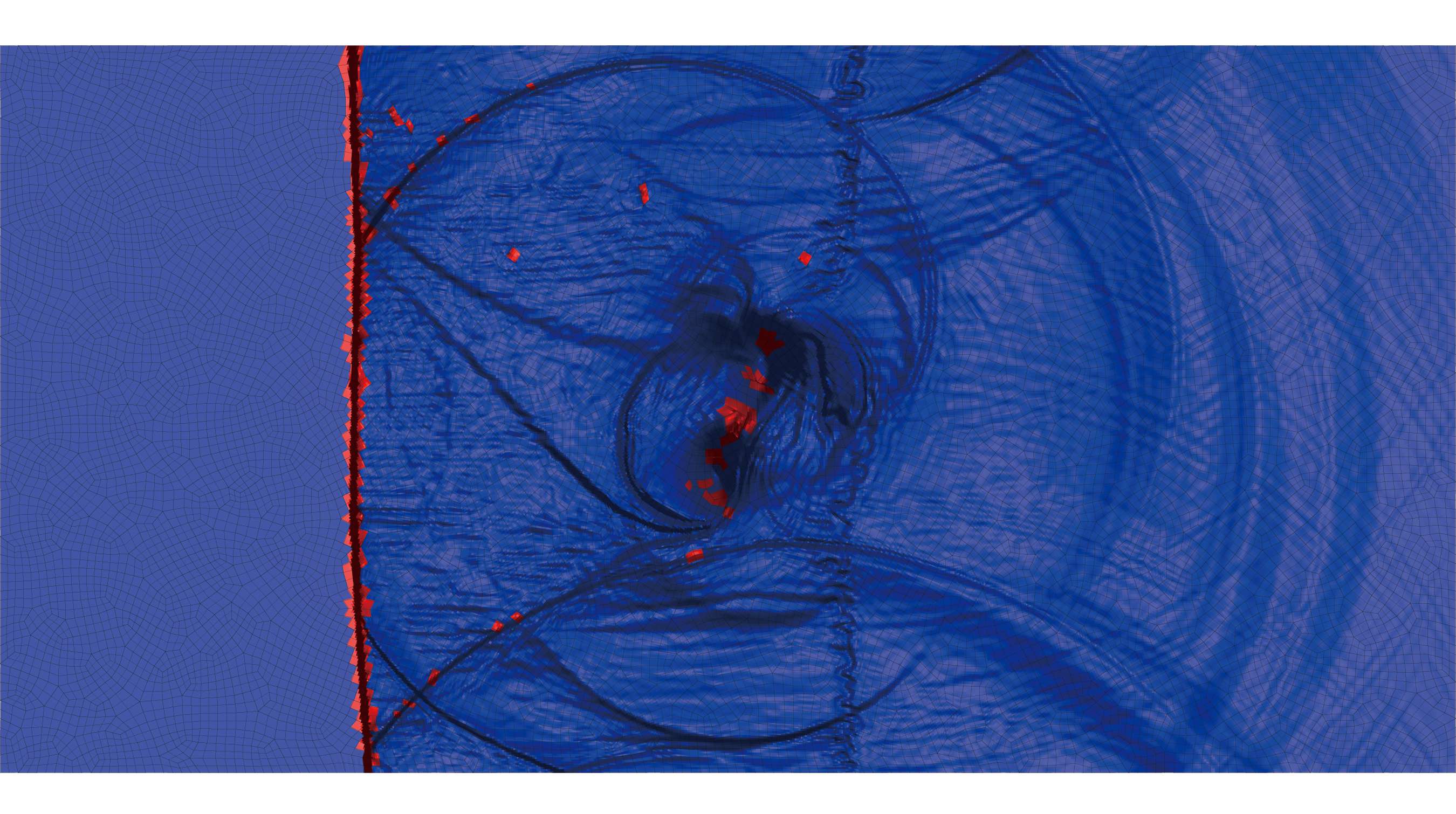}}
		\caption{Results calculated by the hybrid CPR-CNNW2 scheme. T = 7.0. (a)Density. (b)The distribution of troubled cells: the red area represents troubled cells, and the blue area represents smooth cells.}		
		\label{strongshockvortex}
	\end{figure}
	\subsection{Strong shock–vortex interaction problem}
	And the last test, we consider the interaction of traveling vortex with a steady shock proposed by \cite{Rault2003}, unlike section \ref{The ShockVortex Interaction Problem}, the vortex in this section is composite. 
	The computational domain is $\Omega = [0,2]\times[0,1]$ and the initial conditions are given by a stationary shock Mach $M_s$ located at (0.5,y) and by an isentropic vortex centers at $(x_c, y_c) = (0.25,0.5)$. The shock Mach number is denoted by $M_s$ and the strength of the vortex is described in terms of $M_v = v_m/c_0 = v_m/\sqrt{\gamma}$.  We select $M_v = 0.9$ for the vortex and $M_s = 1.5$ for the shock in \cite{Rault2003}. 
	We consider inflow boundary conditions at the left boundary and outflow at the left boundary and the slip-wall boundary conditions for the other boundaries.
\par
	The results on unstructured meshes (cells: 16372) are shown in Figure \ref{strongshockvortex}. Figure \ref{strongshockvortex-rho} and Figure \ref{strongshockvortex-beta} show the distribution of the density and the troubled cells, respectively. These results confirm the ability of the hybrid scheme in capturing shock waves as well as capturing smooth vortex features at the same time. The calculations demonstrate the good properties of the hybrid scheme and are comparable to those in \cite{Rault2003} and \cite{Dumbser2014}.
	
\section{Concluding remarks}\label{concluding} 
	In this paper, we present a priori subcell limiting strategy for the CPR method on unstructured quadrilateral meshes. 
	An indicator based on modal energy coefficient decay is modified and numerical tests show that the modified indicator can detect the discontinuity appearing in cell interfaces. 
	The subcell limiting uses a subcell decomposition based on the nonuniform solution points of CPR for troubled cells. A finite-difference shock-capturing scheme based on nonuniform nonlinear weighted interpolation is developed for the troubled cells. 
	CPR with a priori subcell limiting is a hybrid scheme.
	Numerical results show that the hybrid scheme can capture the strong shock and has a higher resolution than pure second-order CNNW2. The proportion of the troubled cells in the hybrid scheme is small which makes the scheme efficient and has higher resolution.
\section*{Acknowledgement}
	This study was supported by the Basic Research Foundation of National Numerical Wind Tunnel Project, the National Natural Science Foundation of China (Grant Nos. 11902344), and the foundation of State Key Laboratory of Aerodynamics (Grant No. SKLA2019010101). 

\bibliographystyle{unsrt}
\bibliography{Refs}
\end{document}